\RequirePackage{ifpdf}
\ifpdf 
\documentclass[pdftex]{sigma}
\else
\documentclass{sigma}
\fi

\usepackage{amssymb}
\usepackage{paralist}
\usepackage{graphicx}
\usepackage[all]{xy}

\newtheorem*{thma}{Theorem A}
\newtheorem*{thmb}{Theorem B}

\def\mb{\mathbf}

\def\mr{\mathrm}
\def\mc{\mathcal}

\newlength{\equwidth}
\def\A{\mathcal A}
\def\D{\mathcal D}

\def\X{\mathfrak X}
\def\al{\alpha}

\def\ze{\zeta}

\def\ka{\kappa}

\def\la{\lambda}
\def\rh{\rho}

\def\si{\sigma}

\def\Ups{\Upsilon}
\def\ph{\varphi}

\def\om{\omega}
\def\Ga{\Gamma}

\def\Th{\Theta}
\def\La{\Lambda}

\def\Om{\Omega}

\def\pa{\partial}
\def\t{\otimes}

\def\del{\pa}
\def\delstar{\pa^*}

\def\goesto{\rightarrow}

\def\equiv{\Leftrightarrow}
\def\embed{\hookrightarrow}

\def\cinf{\ensuremath{\mathrm{C}^\infty}}

\def\na{\mathrm{\nabla}}

\def\lapl{\mathrm{\square}}
\def\triang{\mathrm{\bigtriangleup}}
\def\Lapl{\mathrm{\triang}}

\def\ker{\mathrm{ker}\ }
\def\im{\mathrm{im}\ }
\def\dim{\mathrm{dim}}
\def\G{\mathcal{G}}
\def\P{\mathcal{P}}
\def\D{\mathcal{D}}
\def\id{\mathrm{id}}

\def\gr{\mathrm{gr}}
\def\id{\mathrm{id}}

\def\Ad{\mathrm{Ad}}

\def\GL{\mathrm{Gl}}
\def\CO{\mathrm{CO}}

\def\ddtz{{\frac{d}{dt}}_{|t=0}}

\def\ti{\tilde}

\def\nn{\ensuremath{\mathbb{N}}}
\def\zz{\ensuremath{\mathbb{Z}}}

\def\rr{\ensuremath{\mathbb{R}}}

\def\CC{\ensuremath{\mathcal{C}}}
\def\ZZ{\ensuremath{\mathcal{Z}}}
\def\BB{\ensuremath{\mathcal{B}}}
\def\HH{\ensuremath{\mathcal{H}}}

\def\VV{\ensuremath{\mathcal{V}}}
\def\TT{\ensuremath{\mathcal{T}}}

\def\SO{\ensuremath{\mathrm{SO}}}

\def\O{\ensuremath{\mathrm{O}}}
\def\GL{\ensuremath{\mathrm{GL}}}

\def\co{\ensuremath{\mathfrak{co}}}

\def\so{\ensuremath{\mathfrak{so}}}
\def\su{\ensuremath{\mathfrak{su}}}

\def\g{\ensuremath{\mathfrak{g}}}

\def\p{\ensuremath{\mathfrak{p}}}
\def\gl{\ensuremath{\mathfrak{gl}}}

\def\ce{\ensuremath{\mathcal{E}}}
\newcommand{\bg}{\mbox{\boldmath{$ g$}}}

\def\Hol{\mathrm{Hol}}

\def\D{\mc{D}}
\def\ti{\tilde}

\begin{document}


\renewcommand{\thefootnote}{$\star$}

\renewcommand{\PaperNumber}{081}

\FirstPageHeading

\ShortArticleName{Conformal Structures Associated to Generic 2-Distributions}

\ArticleName{Conformal Structures Associated to Generic Rank 2\\ Distributions on 5-Manifolds -- Characterization\\ and Killing-Field Decomposition\footnote{This paper is a
contribution to the Special Issue ``\'Elie Cartan and Dif\/ferential Geometry''. The
full collection is available at
\href{http://www.emis.de/journals/SIGMA/Cartan.html}{http://www.emis.de/journals/SIGMA/Cartan.html}}}

\Author{Matthias HAMMERL and Katja SAGERSCHNIG}

\AuthorNameForHeading{M.~Hammerl and K.~Sagerschnig}

\Address{Faculty of Mathematics, University of Vienna, Nordbergstrasse 15, 1090 Vienna, Austria}

\Email{\href{mailto:matthias.hammerl@univie.ac.at}{matthias.hammerl@univie.ac.at}, \href{mailto:katja.sagerschnig@univie.ac.at}{katja.sagerschnig@univie.ac.at}}

\ArticleDates{Received April 09, 2009, in f\/inal form July 28, 2009;  Published online August 04, 2009}

\Abstract{Given a maximally non-integrable $2$-distribution $\D$ on a $5$-manifold   $M$, it was discovered by P. Nurowski that one can naturally associate a
  conformal structure $[g]_{\D}$ of signature $(2,3)$ on $M$. We
  show that those conformal structures $[g]_{\D}$ which come about by
  this construction are characterized by the existence of a
  normal conformal Killing 2-form which is locally decomposable and
  satisf\/ies a genericity condition. We further show that every
  conformal Killing f\/ield of $[g]_{\D}$ can be decomposed into
  a symmetry of $\D$ and an almost Einstein scale of $[g]_{\D}$.}

\Keywords{generic distributions; conformal geometry; tractor calculus; Fef\/ferman construction; conformal Killing f\/ields; almost Einstein scales}

\Classification{34A26; 35N10; 53A30; 53B15; 53B30}

\begin{flushright}
\begin{minipage}{12cm}
\it Dedicated to Peter Michor on the occasion of his 60th birthday
celebrated at the Central European Seminar in Mikulov, Czech Republic,
May 2009.
\end{minipage}
\end{flushright}


\section{Introduction and statement of results}

In this section we brief\/ly introduce the main objects of
interest and state the results of this text.

\subsection{Generic rank 2-distributions on 5-manifolds}\label{sub-gen2}

Let $M$ be a smooth $5$-dimensional manifold and consider
a subbundle~$\D$ of the tangent bund\-le~$TM$ which shall be
of constant rank $2$. We  say that $\D$ is \emph{generic}
if it is maximally non-integrable in the following sense:
For two subbundles $\D_1\subset TM$ and $\D_2\subset TM$
we def\/ine
\begin{gather}\label{distspan}
  [\D_1,\D_2]_x:=\mr{span}(\{[\xi,\eta]_x :\ \xi\in\Ga(\D_1),\eta\in\Ga(\D_2)\}).
\end{gather}
Then we demand that $[\D,\D]\subset TM$  is a subbundle of constant rank $3$  and $[\D,[\D,\D]]=TM$. In other words,  two steps of taking Lie brackets of sections of $\D$ yield all of $TM$.

It is a classical result of \'{E}lie Cartan \cite{cartan-cinq}
that generic rank $2$-distributions on $M$ can equivalently be
described as parabolic geometries of type $(G_2,P)$. This will
be explained in Section~\ref{rank2dist}.

\subsection{The associated conformal structure of signature (2,3)\\ and its characterization}\label{sub-char}

In \cite{nurowski-metric} P.~Nurowski used Cartan's description of generic
distributions to associate to every such distribution a conformal class
$[g]_{\D}$ of signature $(2,3)$-metrics on $M$.

There is a  well studied result similar to Nurowski's construction:  This is
the classical Fef\/fer\-man construction  \cite{fefferman,bds-distinguished,cap-cr-fefferman,cap-gover-holonomy-cr,cap-gover-cr-tractors,felipe-aboutcomplex,felipe-unitaryholonomy}  of a (pseudo) conformal structure on an $S^1$-bundle
over a $CR$ manifold.  It has been observed by A. \v{C}ap in \cite{cap-constructions} that both Nurowski's and Fef\/ferman's results admit interpretations as special cases of a more general construction relating parabolic geometries of dif\/ferent types. In Section~\ref{fefferman}  we will discuss  conformal structures associated to rank two distributions in this picture.
Furthermore, we prove that given a
holonomy reduction of a conformal structure $[g]$ of signature $(2,3)$ to a subgroup of $G_2$, the conformal class $[g]$ is induced by a distribution $\D\subset TM$.

Using strong techniques from the BGG-machinery \cite{BGG-2001,BGG-Calderbank-Diemer,mrh-bgg} and tractor bundles \cite{cap-gover-tractor,cap-gover-irr_tractor},
we then proceed to prove our f\/irst main result in Section~\ref{section-characterization}.
Before we can state it we introduce some simple notation for
tensorial expressions; this is slightly redundant
since we will later use a form of index notation for such formulas.
Take some $g\in [g]$ and
let $\eta\in \otimes^k T^*M$ for $k\geq 2$.
The trace over the $i$-th and the $j$-th slot
of $\eta$ via the (inverse of) the metric $g$
will be denoted $\mr{tr}_{i,j}(\eta)\in\otimes^{k-2} T^*M$.
For an arbitrary tensor $\eta$, $\mr{alt}(\eta)$ is
the full alternation of $\eta$.

\begin{thma}
  Let $[g]$ be a conformal class of signature $(2,3)$ metrics on $M$.
  Then $[g]$ is induced from a generic rank $2$ distribution $\D\subset TM$
  if and only if there exists a normal conformal Killing $2$-form $\phi$ that
  is locally decomposable and satisfies the following genericity condition:
  Let $g\in[g]$ be a metric in the conformal class, $D$ its Levi-Civita connection
  and $P$ its Schouten tensor \eqref{def-schouten}.
  Define
  \begin{gather*}
    \mu :=\mr{tr}_{1,2}{D}\phi\in T^*M, \\
    \rh := +2\Lapl\phi
          +4 \mr{alt}(\mr{tr}_{1,3}DD\phi)
          +3 \mr{alt}(\mr{tr}_{2,3}DD\phi)\\
          \phantom{\rh :=}{}
          +24 \mr{alt}(\mr{tr}_{1,3}P\t\phi)
          -6 \mr{tr}_{1,2}P\t\phi\in \La^2 T^*M.
  \end{gather*}
  Here we use the convention
  \begin{gather*}
    \Lapl\si=-\mr{tr}_{1,2}DD\si.
  \end{gather*}
  Then one must have
  \begin{gather*}
    \phi\wedge\mu\wedge\rh\not=0.
  \end{gather*}
\end{thma}

To be precise, Theorem~A assumes orientability of $TM$, but
this is a minor assumption only made for convenience
of presentation -- see Remark~\ref{orientability}.

\subsection{Killing f\/ield decomposition}\label{sub-kill}

Let $\mb{sym}(\D)$ denote the inf\/initesimal symmetries of the distribution
$\D$, i.e.,
\begin{gather*}
  \mb{sym}(\D)=\{\xi\in\X(M) :\ \mc{L}_{\xi}\eta=[\xi,\eta]\in\Ga(\D)\ \forall\, \eta\in\Ga(\D)\}.
\end{gather*}
The corresponding objects for conformal structures
are the \emph{conformal\ Killing\ fields}
\begin{gather*}
  \mb{cKf}([g])=\{\xi\in\X(M)  :\ \mc{L}_{\xi}g=e^{2f}g\ \mr{for\ some}\ g\in[g]\ \mr{and}\ f\in\cinf(M)\}.
\end{gather*}
Since the construction that associates a conformal structure $[g]_{\D}$ to a distribution $\D$ is natural, every
symmetry $\xi\in\X(M)$ of the distribution will also preserve
the associated conformal structure $[g]_{\D}$, i.e., it is  a conformal Killing f\/ield. This yields an embedding
 \begin{gather*}
   \mb{sym}(\D)\embed\mb{cKf}([g]_{\D}).
 \end{gather*}
We will show  that a complement to
$\mb{sym}(\D)$ in $\mb{cKf}([g]_{\D})$ is given by the \emph{almost\ Einstein\ scales}\   of $[g]_{\D}$: a function $\si\in\cinf(M)$ is
an almost Einstein scale
for $g\in[g]_{\D}$ if it is non-vanishing on an open dense subset $U$ of $M$\
and satisf\/ies that $\si^{-2}g$ is Einstein on $U$~\cite{gover-aes}. The natural origin
of almost Einstein scales via tractor calculus will be seen in
Section~\ref{section-conformal}.

\begin{thmb} Let $[g]_{\D}$ be the conformal structure associated to a generic rank~$2$ distribution~$\D$ on a $5$-manifold $M$, and let $\phi$ be a conformal Killing form characterizing the conformal structure as in Theorem~A.
      Then every conformal Killing field decomposes into a symmetry of
    the distribution~$\D$ and an almost Einstein scale:
    \begin{gather}
    \label{decomp1}
      \mb{cKf}([g]_{\D})=\mb{sym}(\D)\oplus \mb{aEs}([g]_{\D}).
    \end{gather}

    The mapping that associates to a conformal Killing field $\xi\in\X(M)$
    its
    almost Einstein scale part with respect to the decomposition \eqref{decomp1}
    is given by
\begin{gather*}
      \xi\mapsto \mr{tr}_{1,3}\mr{tr}_{2,4}\left(\phi\t D\xi
-\tfrac{1}{2}\xi\t D\phi\right),
    \end{gather*}
    where $D$ is the Levi-Civita connection of an arbitrary metric
    $g$ in the conformal class.

    The mapping that associates to an almost Einstein scale $\si\in\cinf(M)$
    $($for a metric $g\in[g])$ a conformal Killing field is given by
\begin{gather*}
      \si\mapsto \mr{tr}_{2,3}\phi\t (D\si)-\tfrac{1}{4}\mr{tr}_{1,2}(D\phi)\si.
    \end{gather*}
\end{thmb}

We remark here that the constructions in this paper, both for characterization (Section~\ref{section-characterization}) and automorphism-decomposition (Section~\ref{section-decomp}),
are largely analogous to the ones of \cite{cap-gover-holonomy-cr}
and~\cite{cap-gover-cr-tractors} for the (classical) Fef\/ferman spaces.

This article incorporates material of the authors' respective theses~\cite{katja-thesis,mrh-thesis}.

\section{Preliminaries on Cartan and parabolic geometries}\label{section-general}
In this section we discuss general parabolic geometries.
These are special kinds of Cartan geo\-met\-ries.

\subsection{Cartan geometries}
  Let $G$ be a Lie group and $P<G$ a closed subgroup.
  The Lie algebras of $P$ and $G$ will be denoted by $\p$ and $\g$.
  Let $\G\overset{\pi}{\goesto}M$ be a $P$-principal bundle over a manifold $M$.
  The right action of $P$ on $\G$ will be denoted by
 $r^p(u)=u\cdot p$\
    for $u\in\G$ and $p\in P$. The corresponding
    fundamental vector f\/ields are $\ze_Y(u):=\ddtz r^{\exp(tY)}(u)=\ddtz u\cdot \exp(tY)$
   for $Y\in\p$.
\begin{definition}
  A \emph{Cartan geometry} of type $(G,P)$ on a manifold $M$ is
  a $P$-principal bundle $\G\overset{\pi}{\goesto}M$ endowed
  with a \emph{Cartan connection form}\ $\om\in\Om^1(\G,\g)$,
  i.e., a $\g$-valued $1$-form on $\G$ satisfying
  \begin{enumerate}[(C.1)]\itemsep=0pt
  \item\label{cartan1} $\om_{u\cdot p}(T_u r^p \xi)=\Ad(p^{-1})\om_u(\xi)$ for all $p\in P$, $u\in \G,$ and $\xi\in T_u\G$.
  \item\label{cartan2} $\om(\ze_Y)=Y$ for all $Y\in\p$.
  \item\label{cartan3} $\om_u:T_u\G\goesto\g$ is an isomorphism for all $u\in\G$.
  \end{enumerate}
 We say that $\om$ is right-equivariant, reproduces fundamental vector f\/ields
  and is an absolute parallelism $\om:T\G\cong \G\times \g$.
\end{definition}

It is easily seen that for a Cartan geometry $(\G,\om)$ of
type $(G,P)$ the map $\G\times\g\mapsto TM$ given by $(u,X)\mapsto T_u \pi \omega_u^{-1}(X)$ induces an
isomorphism
\begin{gather*}
  TM=\G\times_P\g/\p.
\end{gather*}
In particular, $\dim\ M=\dim\ \g/\p$.

Cartan geometries can be viewed as curved versions of homogeneous spaces: The \emph{homogeneous model} of Cartan geometries of type $(G,P)$ is the  principal bundle $G\goesto G/P$ endowed with the Maurer--Cartan form $\om^{MC}\in\Omega^1(G,\g)$, which satisf\/ies the Maurer--Cartan equation \[d\om^{MC}(\xi,\eta)+\big[\om^{MC}(\xi),\om^{MC}(\eta)\big]=0\] for all $\xi,\eta\in\X(G)$.

For a general Cartan geometry $(\G,\om)$  the failure of $\om$ to satisfy the Maurer--Cartan equation
is measured by the  \emph{curvature}\ form $\Om\in\Om^2(\G,\g)$,
\begin{gather}\label{def-curvature}
   \Om(\xi,\eta)=d\om(\xi,\eta)+[\om(\xi),\om(\eta)].
\end{gather}
One can show that $\Om$ vanishes, i.e., $\om$ is f\/lat, if and only if
$(\G,\om)$ is locally isomorphic (in the obvious sense) with $(G,\om^{MC})$.

Since the Cartan connection def\/ines an absolute parallelism $\om:T\G\cong \G\times \g$, its curvature can be equivalently encoded in the curvature function $\ka\in\cinf(\G,\La^2(\g^*)\t\g)$, \[\ka(u)(X,Y):=\Om\big(\om_u^{-1}(X),\om_u^{-1}(Y)\big).\] One  verif\/ies that $\Om$  vanishes on vertical f\/ields $\ze_Y$ for $Y\in\p$, i.e., it is horizontal. This implies, that $\kappa$ in fact def\/ines  a function $\mathcal{G}\mapsto\Lambda^2(\g/\p)^*\otimes\g$. And since $\Om$ is $P$-equivariant, so
is $\ka$.

We denote by $\A M:=\G\times_P\g$
the associated bundle corresponding to the restriction of the adjoint representation $\Ad:G\goesto\GL(\g)$ to $P$. It is  called the \emph{adjoint tractor bundle} (general tractor bundles will be
introduced below).

Note that since the curvature of a Cartan connection is horizontal and $P$-equivariant, it factorizes to a $\A M$-valued 2-form $K\in\Om^2(M,\A M)$ on $M$.
Thus, $\Om\in\Om^2(\G,\g),K\in\Om^2(M,\A M)$ and $\ka:\G\goesto\La^2(\g/\p)^*\t\g$ all encode essentially the same object, namely the curvature of the
Cartan connection form $\om$, and technical reasons will determine
which representation should be used at a given point.

\subsection{Tractor bundles}\label{tractors}
For any $G$-representation $V$, the associated bundle
\begin{gather*}
  \VV:=\G\times_P V
\end{gather*}
 is called a \emph{tractor\ bundle}.
 Tractor bundles  carry canonical linear connections:
Extend the structure group of $\G$ from $P$ to $G$\
by forming $\G':=\G\times_P G$. Then $\om$ extends uniquely
to a~$G$-equivariant $\g$-valued $1$-form $\om'$ on $\G$ reproducing fundamental
vector f\/ields -- i.e.,  to  a~principal connection form.
Since $\VV=\G\times_P V=\G'\times_G V$ one has
the induced \emph{tractor\ connection}\ $\na^{\VV}$ on~$\VV$.
For computations we will use the following explicit formula:
Let $\xi'\in\X(\G)$ be a $P$-invariant lift of a vector f\/ield $\xi\in\X(M)$\
and $f_s\in\cinf(\G,V)^P$ be the $P$-equivariant $V$-valued function
on $\G$ corresponding to $s\in\Ga(\VV)$. Then $\na_{\xi}^{\VV}s$\
corresponds to the $P$-equivariant function
\begin{gather}\label{deftracon}
  \xi'\cdot f_s+\om(\xi') f_s.
\end{gather}

\subsection{Parabolic geometries}

We now specialize to parabolic geometries.
This class of Cartan geometries has a natural algeb\-raic normalization condition
which is employed in this text to describe conformal structures and generic
distributions as parabolic geometries. For a thorough discussion
of parabolic geometries we refer to \cite{cap-slovak-par}.

Let us start with the algebraic background and introduce the notion of a  $|k|$-graded Lie algebra~$\g$:  this is a a semisimple Lie algebra together with a vector space decomposition $\g=\g_{-k}\oplus\cdots\oplus \g_k$  such that $[\g_i,\g_j]\subset \g_{i+j}$.
Then $\g_-=\g_{-k}\oplus\cdots\oplus \g_{-1}$\
and $\p_+=\g_1\oplus\cdots\oplus\g_k$ are nilpotent subalgebras of $\g$. The Lie algebra $\p=\g_0\oplus\cdots\oplus\g_k$ is indeed a parabolic subalgebra of~$\g$, and~$\g_0$ is its reductive Levi part.  The grading induces a f\/iltration on $\g$ via
$\g^i:=\g_i\oplus\cdots\oplus\g_k$.

Let $G$ be a Lie group with Lie algebra $\g$.
Let $P$  be a closed subgroup of $G$ whose Lie algebra is the parabolic $\p\subset \g$ and such that it preserves the f\/iltration, i.e., for all $p\in P$  we have $\Ad(p)\g^i\subset\g^i\ \forall i \in\zz$.

\begin{definition}
 A Cartan geometry of type $(G,P)$ for groups as introduced above is called a~\emph{parabolic geometry}.
\end{definition}
In the following, we will also consider the subgroup
\begin{gather*}
  G_0:=\{g\in P: \ \Ad(g)\g_i\subset \g_i\ \forall\, i\in \zz\}.
\end{gather*}
of all elements in $P$ preserving the grading on the Lie algebra, which has Lie algebra $\g_{0}$, and the subgroup
\[P_+:=\{p\in P: \ (\Ad(p)-\id)\g^i\subset\g^{i+1}\ \forall\, i\in\zz\},\] which has Lie algebra $\p_{+}$.
Actually $P$ decomposes as a semidirect product
\begin{gather*}
  P=G_0\ltimes P_+;
\end{gather*}
thus $P/P_+=G_0$.

\subsection{Lie algebra dif\/ferentials and normality}

Let $V$ be a $G$-representation.
We now introduce algebraic dif\/ferentials \[\del:\La^i(\g/\p)^*\t V\goesto
\La^{i+1}(\g/\p)^*\t V\] and \[\delstar:\La^{i+1}(\g/\p)^*\t V\goesto
\La^{i}(\g/\p)^*\t V.\] For the f\/irst of these, we use the  $G_0$-equivariant
identif\/ication of $\g_-^*$ with $(\g/\p)^*$ and def\/ine $\partial$ as the dif\/ferential computing the Lie algebra cohomology of $\g_{-}$ with values in $V$.
For the  \emph{Kostant codifferential} $\partial^*$ we use the $P$-equivariant identif\/ication of $(\g/\p)^*$ with $\p_+$ given by the Killing form; it is then def\/ined as the dif\/ferential computing the Lie algebra homology of $\mathfrak{p}_{+}$ with values in $V$.
We include the explicit formula for $\partial^*$, which will be needed later on:
On a decomposable element $ \ph=Z_1\wedge\cdots\wedge Z_i \t v\in\La^i\p_+\t V,$   $\partial^*$ is given as
\begin{gather*}
   \delstar(\ph) :=\sum\limits_{j=1}^i (-1)^j Z_1\wedge \cdots\wedge \widehat {Z_j}\wedge\cdots \wedge Z_i \t(X_j v)\\
\phantom{\delstar(\ph) :=}{} +\sum\limits_{1\leq j<k\leq i}(-1)^{j+k} [Z_j,Z_k]\wedge Z_1\wedge \cdots\wedge\widehat{Z_j}\wedge\cdots \wedge\widehat{Z_k}\wedge\cdots \wedge Z_i\otimes v.
\end{gather*}

The $G$-representation $V$ carries a natural
$G_0$-invariant grading $V_0\oplus\cdots \oplus V_r$ for some $r\in \nn$.
The induced f\/iltration $V^i:=V_i\oplus\cdots \oplus V_r$ is even $P$-invariant.
The grading on $V$ and the grading $\g_{-k}\oplus\cdots\oplus\g_{-1}=\g_-\cong\g/\p$ naturally induce a grading on the spaces $C_i(V)=\La^i (\g/\p)^*\t V$, which is preserved by both $\del$ and $\delstar$. While $\delstar$ is seen to be $P$-equivariant, $\del$ is only $G_0$-equivariant.

We consider the spaces of cocycles $Z_i(V):=\ker\delstar\subset \La^i(\g/\p)^*\otimes V$, coboundaries $B_i(V):=\im\delstar\subset \La^i(\g/\p)^*\otimes V$ and cohomologies $H_i(V):=Z_i(V)/B_i(V)$. Let $\Pi_i:Z_i(V)\to H_i(V)$ be the canonical surjection. By~\cite{kostant-61}, the dif\/ferentials $\partial$ and $\partial^*$ are adjoint with respect to a natural inner product on the space  $C_i(V)$. Via the Kostant Laplacian
$\lapl=\delstar\circ\del+\del\circ\delstar$
this yields a~$G_0$-invariant Hodge decomposition
\begin{gather*}
  C_i(V)=\im\del\oplus\ker\lapl\oplus\im\delstar.
\end{gather*}
Thus, as a $G_0$-module, $H_i(V)$ can be embedded into $Z_i(V)\subset C_i(V)$.

Since $\G\times_P\g/\p=T M$ the associated bundle $\CC_i:=\G\times_P\La^i(\g/\p)^*\t  V$\
of $C_i(V)$ is $\La^i T^*M \t \VV$, whose sections are the $\VV$-valued
$i$-forms $\Om^i(M,\VV)$. The $P$-equivariant dif\/ferential
$\delstar$ carries over to the associated spaces,
\begin{gather*}
  \delstar:\Om^{i+1}(M,\VV)\goesto \Om^i(M,\VV).
\end{gather*}
We set $\ZZ_i(\VV):=\G\times_P Z_i(V)$, $\BB_i(\VV):=\G\times_P B_i(V)$
and $\HH_i(\VV):=\G\times_P H_i(V)$. The canonical
surjection from $\ZZ_i(\VV)$ onto $\HH_i(\VV)$ is denoted by $\Pi_i$.
If the tractor bundle $\VV$ in question is unambiguous we just write
$\CC_i$, $\ZZ_i$, $\BB_i$, $\HH_i$.

The Kostant codif\/ferential provides a conceptual normalization condition for parabolic geometries:
 Recall that the curvature  of a Cartan connection form $\om$ factorizes to a two-form $K\in\Om^2(M,\A M)$ on $M$ with values
in the adjoint tractor bundle.
\begin{definition}\label{def-normal}
 A Cartan connection form $\om$ is called \emph{normal}\ if $\delstar(K)=0$. In this case one has
    the \emph{harmonic\ curvature}\ $K_H=\Pi_2(K)\in\HH_2(\A M)$.
\end{definition}
In the picture of $P$-equivariant functions on $\mathcal{G}$, the harmonic curvature corresponds to the composition of the curvature function $\kappa$ with the projection $\Pi_2:Z_2(\g)\to H_2(\g)$, i.e., to $\kappa_{H}=\Pi_2\circ\kappa$ .

There is a simple algorithm to compute the cohomology spaces $H_i(V)$ provided
by Kostant's version of the Bott--Borel--Weil theorem, cf.~\cite{kostant-61,silhan-bott}.
Mostly, we will just need to know $H_0(V)$, which turns out
to be  the \emph{lowest\ homogeneity} of $V$, i.e.,
$H_0(V)=V/V^1=V/(\p_+V)$.

\subsection{The BGG-(splitting-)operators}\label{section-bgg}

The BGG-machinery developed
in \cite{BGG-2001}\ and \cite{BGG-Calderbank-Diemer}\ will
feature prominently at many crucial points in this paper.
The presentation here is very brief and the most
important operators will later be given explicitly (see the end
of the next section on conformal geometry). The
highly useful Lemma \ref{reproduceparlem}\ below can be understood
without its relation to the BGG-machinery.

The main observation  is that
for every $\si\in\Ga(\HH_0)$ there is a unique $s\in \Ga(\VV)$ with
$\Pi_0(s)=\si$\
such that $\na^{\VV}s\in\Ga(\ZZ_{1})$, i.e., such that $\delstar(\na^{\VV}s)=0$.
This gives a natural splitting
$L_0^{\VV}:\Ga(\HH_0)\goesto\Ga(\VV)$\
of $\Pi_0:\Ga(\VV)\goesto\Ga(\HH_0)$ called the $1$-st\ \emph{BGG-splitting operator} and it def\/ines the
$1$-st \emph{BGG-operator}
\begin{gather*}
  \Th_0^{\VV} : \ \Ga(\HH_0)\goesto\Ga(\HH_{1}),\qquad
  \si \mapsto\Pi_{1}(\na^{\VV}s(L_0(\si))).
\end{gather*}
We only remark that this construction of dif\/ferential splitting operators of
the projections $\Pi_i:\ZZ_i\goesto\HH_i$ proceeds similarly,
and one obtains the celebrated BGG-sequence $\HH_i\overset{\Th_i}{\goesto}\HH_{i+1}$.

We will often need the following consequence of
the def\/inition of $L_0^{\VV}$:
If $s\in\Ga(\VV)$ is parallel, one trivially has $\delstar(\na^{\VV}s)=0$,
and thus $s=L_0^{\VV}(\Pi_0(s))$. This is important enough
to merit a
\begin{lemma}\label{reproduceparlem}
  On the space of parallel sections of a tractor bundle $\VV$,
  $L_0^{\VV}\circ \Pi_0$ is the identity, i.e.,
  if $s\in\Ga(\VV)$ with $\na^{\VV}s=0$, then
\begin{gather*}
  s=L_0^{\VV}(\Pi_0(s)).
\end{gather*}
In particular, if the projection of a parallel section $s\in\Ga(\VV)$ to its  part in $\Ga(\HH_0)=\Ga(\VV/\VV^1)$ vanishes,
$s$ must already have been trivial.
\end{lemma}

We now proceed to discuss conformal structures as parabolic geometries
in Section~\ref{section-conformal}\ and do likewise for generic rank
two distributions in Section~\ref{rank2dist}.

\section{Conformal structures}\label{section-conformal}

Two pseudo-Riemannian metrics $g$ and $\hat g$ with
signature $(p,q)$ on a $n=p+q$-dimensional manifold
$M$ are said to be \emph{conformally equivalent} if there
is a function $f\in\cinf(M)$ such that $\hat g=e^{2f}g$.
The conformal equivalence class of $g$ is denoted by $[g]$
and $(M,[g])$ is said to be a~manifold endowed with a conformal structure.
An equivalent description of a conformal
structure of signature $(p,q)$ is a a reduction of structure group
of $TM$ to $\CO(p,q)=\rr_+\times \O(p,q)$,
and the corresponding $\CO(p,q)$-bundle will be denoted ${\ti\G}_0$.

The associated bundle to ${\ti\G}_0$ for the $1$-dimensional
representation $\rr[w]$ of $\CO(p,q)$ given by
\begin{gather*}
  (c,C)\in \CO(p,q)=\rr_+\times O(p,q)\mapsto c^{w}
\end{gather*}
for $w\in\rr$ is called the bundle of conformal $w$-densities
and denoted by $\ce[w]$.

We will use abstract index notation and notation for weighted bundles
similar to \cite{gover-silhan-2006}: $\ce_a:=T^*M$, $\ce^a:=TM$,
$\ce_{ab}=\ce_a\t\ce_b$, $\ce_a[w]:=\ce_a\t\ce[w]$.
Recall the Einstein convention, e.g., for $\xi^a\in\Ga(\ce^a)=\X(M)$ and $\ph_a\in\Ga(\ce_a)=\Om^1(M)$, $\xi^a\ph_a=\ph(\xi)\in\cinf(M)$.
Round brackets will denote symmetrizations, e.g.~$\ce_{(ab)}=S^2T^*M$
and square brackets anti-symmetrizations, e.g.~$\ce_{[ab]}=\La^2 T^*M$.
In the following we will not distinguish between the
space of sections $\Ga(\ce_{a\cdots b}[w])$ and $\ce_{a\cdots b}[w]$ itself.

Given a metric $g\in [g]$, a section $\si\in\ce[w]$\
trivializes to a function $[\si]_g\in\cinf(M)$ and one has
\begin{gather*}
  [\si]_{e^2f g}=e^{wf}[\si]_g.
\end{gather*}
Tautologically, the conformal class of metrics $[g]$ def\/ines
a canonical section $\bg$ in $\ce_{(ab)}[2]=\Ga(S^2T^*M\t\ce[2])$, called the \emph{conformal\ metric}, such that the
trivialization of $\bg$ with respect to $g\in [g]$ is just $g$.
The conformal metric $\bg$ allows one to raise or lower indices
with simultaneous adjustment of the conformal weight: e.g.,
for a vector f\/ield $\xi^p\in\ce^p=\X(M)$ one can form $\xi_p=\bg_{pq}\xi^q\in\ce_p[2]=\Ga(T^*M\t\ce[2])$, which is a $1$-form of weight $2$.
\subsection{Conformal structures as parabolic geometries}\label{conform-parabolic}
Let $M_{p,q}$ be a given symmetric bilinear form of signature $(p,q)$ on $\rr^n=\rr^{p,q}$, and def\/ine the symmetric bilinear form $h$ of signature $(p+1,q+1)$ on $\rr^{n+2}$ by
\begin{gather}\label{metrich}
  h=
  \begin{pmatrix}
    0 & 0 & 1 \\
    0 & M_{p,q} & 0 \\
    1 & 0 & 0
  \end{pmatrix}.
\end{gather}
We def\/ine  ${\ti P}\subset\SO(h)\cong\SO(p+1,q+1)$ as
the stabilizer of the isotropic ray $\rr_+ e_1$, and one f\/inds
 $\ti P=\CO(p,q)\ltimes {\rr^n}^*$.
The Lie algebra $\so(p+1,q+1)=\so(h)$ is $|1|$-graded
\begin{gather*}
  \so(h)={\so(h)}_{-1}\oplus{\so(h)}_0\oplus {\so(h)}_1=\rr^n\oplus\co(p,q)\oplus {\rr^n}^*.
\end{gather*}
Realized in $\gl(n+2)$ it is given by matrices of the form
\begin{gather}\label{solie}
  \begin{pmatrix}
      -\al & \vline & -Z^tM_{p,q} & \vline & 0 \\
  \hline
  X & \vline & A & \vline & Z
  \\
  \hline
  0 & \vline & -X^tM_{p,q} & \vline & \al
  \end{pmatrix}
  , \qquad \al\in\rr, \quad X,Z\in\rr^n, \quad A\in\so(M_{p,q}).
\end{gather}
Let $(\ti\G\goesto M,\ti\om)$
be a Cartan geometry of type $(\SO(h),\ti P)$.
Def\/ine
\begin{gather*}
  \ti\G_0:=\ti\G/\ti P_+=\ti\G/{\rr^n}^*.
\end{gather*}
Then $\ti\G_0$ is a $\CO(p,q)$-principal bundle
over $M$ and \[TM=\ti\G\times_{\ti P}\so(h)/\p=
\ti\G_0\times_{\CO(p,q)}\rr^n,\]
i.e., $\ti\G_0\goesto M$ gives a reduction of
structure group of $TM$ to $\CO(p,q)$ and thus
a conformal structure of signature $(p,q)$.

Since there are
many non-isomorphic Cartan geometries of type $(\SO(h),\ti P)$
describing the same conformal structure on the
underlying manifold, one imposes a normalization
condition on the curvature $K\in\Om^2(M,\ti\A M)$ of $\ti\om$.
Using the notion of normality introduced in Def\/inition~\ref{def-normal},
one has:
\begin{theorem}[\cite{cartan-conformal}]\label{cartanconform}
Up to isomorphism there
is a unique ${\ti P}$-principal bundle ${\ti\G}$ over $M$\
endowed with a normal Cartan connection
form $\ti\om\in\Om^1({\ti\G},\so(h))$ such that ${\ti\G}/{{\rr^n}^*}={\ti\G}_0$
is the conformal frame bundle of $(M,[g])$.
\end{theorem}
This provides an equivalence of categories between
oriented conformal structures of signature $(p,q)$
and normal parabolic geometries of type $(\SO(h),{\ti P})$.

\subsection{Tractor bundles for conformal structures}\label{conftra}

The \emph{standard tractor bundle} of conformal
geometry is obtained by the associated bundle $\TT:={\ti\G}\times_{\ti P} \rr^{n+2}$
of the standard representation of $\ti P=\SO(h)$ on $\rr^{n+2}$.
$\ti P$ preserves the f\/iltration
\begin{gather}\label{rrnfilt}
  \{0\}
  \subset
  \begin{pmatrix}
    \rr\\
    0 \\
    0
  \end{pmatrix}
  \subset
  \begin{pmatrix}
    \rr\\
    \rr^n\\
    0
  \end{pmatrix}
  \subset
  \begin{pmatrix}
    \rr \\
    \rr^n\\
    \rr
  \end{pmatrix}
\end{gather}
of $\rr^{n+2}$, and therefore gives a well-def\/ined f\/iltration
$\{0\}\subset\TT^{1}\subset\TT^0\subset\TT^{-1}=\TT$.
The associated graded of $\TT$ is $\gr(\TT)=\gr_{-1}(\TT)\oplus\gr_0(\TT)\oplus\gr_1(\TT)$, with
\begin{gather}
  \gr_1(\TT):=\TT^1=\ce[-1],\nonumber \\
    \gr_0(\TT):=\TT^0/\TT^1=\ce_a[-1],\label{tragrad} \\
  \gr_{-1}(\TT):=\TT^{-1}/\TT^0=\ce[1].\nonumber
\end{gather}
It is a general and well known fact of conformal tractor calculus that
a choice of metric $g\in [g]$ yields a reduction
of the $\ti P$-principal bundle $\ti\G$ to the
$\CO(p,q)$-principal bundle $\ti\G_0$, and this
is seen to provide an isomorphism
of a natural bundle with its associated graded space.
In the case of the standard tractor bundle $\TT$, this
gives an isomorphism of $\TT$ with $\gr(\TT)$,
and a section $s\in\Ga(\TT)$ will then be written
\begin{gather}\label{std_tractor}
 [ s]_{g}=
  \begin{pmatrix}
    \rh \\
    \ph_a \\
    \si
  \end{pmatrix}
  \in
  \begin{pmatrix}
    \ce[-1] \\
    \ce_a[1] \\
    \ce[1]
  \end{pmatrix}.
\end{gather}
For $\hat g=\mr{e}^{2f}g$ one has the transformation
\begin{gather*}
        [s]_{\hat g}=
        \begin{pmatrix}
            \hat\rh \\
            \hat\ph_a \\
            \hat\si
        \end{pmatrix}=
        \begin{pmatrix}
            \rh-\Ups_a\ph^a-\frac{1}{2}\si\Ups^b\Ups_b \\
            \ph_a+\si \Ups_a  \\
            \si
        \end{pmatrix}
\end{gather*}
where $\Ups=df$.
The insertion of $\ce[-1]$ into $\TT$ as the top slot is independent of
the choice of $g\in[g]$ and def\/ines a section $\tau_+\in\TT[1]$.
The insertion of $\ce[1]$ into $\TT$ as the bottom slot is
well def\/ined only via a choice of $g\in[g]$ and def\/ines a section $\tau_-\in\TT[-1]$. Let $e_1,\dots,e_{n+2}$ be the standard basis of $\rr^{n+2}$.
Then $\tau_+$ and $\tau_-$ can be understood as the sections corresponding
to the constant functions on $\ti\G_0$ mapping to $e_1\t 1\in\rr^{n+2}\t\rr[1]$ resp.\ $e_{n+2}\t 1\in\rr^{n+2}\t\rr[-1]$.

Since $h\in S^2 T^*\rr^{n+2}$ is $\SO(h)$-invariant
it def\/ines a \emph{tractor metric} $\mb{h}$ on $\TT$.
With respect to $g\in[g]$ and the decomposition \eqref{std_tractor}
of an element $s\in\Ga(\TT)$
\begin{gather*}
  [\mb{h}]_g=
  \begin{pmatrix}
    0 & 0 & 1 \\
    0 & \bg & 0 \\
    1 & 0 & 0
  \end{pmatrix}.
\end{gather*}
Let $D$ be the Levi-Civita connection of $g\in [g]$,
then the tractor connection $\na^{\TT}$ on $\TT$ is
given~by
\begin{gather}\label{sttracon}
  [\na^{\TT}_c s]_g=\na^{\TT}_c
  \begin{pmatrix}
    \rh \\
    \ph_a \\
    \si
  \end{pmatrix}
  =
  \begin{pmatrix}
    D_c \rh-P_{c}^{\; b}\ph_b \\
    D_c\ph_a+\si P_{ca}+\rh \bg_{c a}\\
    D_c \si-\ph_c
  \end{pmatrix}.
\end{gather}
Here
\begin{gather}\label{def-schouten}
  P=P(g)=\frac{1}{n-2}\left(\mr{Ric}(g)-\frac{\mr{Sc}(g)}{2(n-1)}g\right)
\end{gather}
is the \emph{Schouten tensor} of $g$.
The trace of
the Schouten tensor is denoted $J=g^{pq}P_{pq}$.

The adjoint tractor bundle is ${\ti\A} M={\ti\G}\times_{\ti P}\so(h)$,
which can be identif\/ied with $\so(\TT,{\mb{h}})=\La^2\TT$. With respect to $g\in [g]$,
${\ti\A} M=TM\oplus \co(TM,g)\oplus T^* M$, and in matrix
notation a section $[s]_g=\xi\oplus (\al,A)\oplus \ph\in\X(M)\oplus\co(TM,g)\oplus \Om^1(M)$ will be written as
\begin{gather*}
  \begin{pmatrix}
      -\al & \vline & -\ph_a & \vline & 0 \\
  \hline
  \xi^a & \vline & A & \vline & \ph^a
  \\
  \hline
  0 & \vline & -\xi_a & \vline & \al
  \end{pmatrix}.
\end{gather*}
The curvature form ${\ti K}\in\Om^2(M,{\ti\A} M)$  has
in fact values in ${\ti\A} M^0$; this is called \emph{torsion-freeness}. It furthermore
decomposes into Weyl curvature $C\in\Om^2(M,\so(TM))=\ce_{c_1c_2\;\; d}^{\;\ \;\; \ c}$ and
Cotton--York tensor $A\in\ce_{a[c_1c_2]}$:
\begin{gather}\label{formulaK}
  {\ti K}_{c_1c_2}=
  \begin{pmatrix}
    0 & -A_{ac_1c_2} & 0 \\
    0 & C_{c_1c_2\;\; b}^{\;\ \;\; \ a} & A^a_{\;\ c_1c_2} \\
    0 & 0 & 0
  \end{pmatrix}.
\end{gather}
The Weyl curvature $C$ is the completely trace-free part of
the Riemannian curvature $R$ of $g$. The Cotton--York tensor
is given by $A=A_{ac_1c_2}=2D_{[c_1}P_{c_2]a}$.

We will later need the f\/irst BGG-splitting operators
for the tractor bundles $\TT,{\ti\A} M=\La^2 \TT$ and $\La^3 \TT$,
and therefore give general formulas from~\cite{mrh-bgg} for the space $\VV:=\La^{k+1}\TT$ for $k\geq 0$.

The $\ti P$-invariant f\/iltration \eqref{rrnfilt} of
$\rr^{n+2}$ from above carries over
to the invariant f\/iltration of the exterior power $V=\La^{k+1}\rr^{n+2}$,
written $\{0\}\subset V^{1}\subset V^0 \subset V^{-1}=V$.
Again, this yields f\/iltrations of the associated bundles:
$\{0\}\subset\VV^{1}\subset\VV^0\subset\VV^{-1}=\VV:=\La^{k+1}\TT$.
The notion of the associated graded space is the same:
we def\/ine $\gr(\La^k\TT)$ as the direct sum over all
$\gr_i(\La^{k}\TT):=(\La^k\TT)^i/(\La^k\TT)^{i+1}$.
With respect to $g\in[g]$, for $k\geq 0$,
one again obtains an isomorphism of $\La^{k+1}\TT$\
with $\gr(\La^{k+1}\TT)$, and we will write
\begin{gather*}
  [\La^{k+1}\TT]_g=
    \begin{pmatrix}
        \ce_{[a_1\cdots a_k]}[k-1] \\
        \ce_{[a_1\cdots a_{k+1}]}[k+1]\; |\;\ \ce_{[a_1\cdots a_{k-1}]}[k-1] \\
        \ce_{[a_1\cdots a_k]}[k+1]
    \end{pmatrix}.
\end{gather*}
This identif\/ication employs the insertions of the top slot $\tau_+\in\TT[1]$ and bottom slot $\tau_-\in\TT[-1]$:
\begin{gather}\label{identification}
      \begin{pmatrix}
            \rh_{a_1\cdots a_k} \\
            \ph_{a_0\cdots a_k}\; |\;\ \mu_{a_2\cdots a_k} \\
            \si_{a_1\cdots a_k}
    \end{pmatrix}
    \mapsto
    \tau_-\wedge\si+\ph+\tau_+\wedge\tau_-\wedge\mu+\tau_+\wedge\rh.
\end{gather}
To understand the map $\si\mapsto \tau_-\wedge\si$ better,
observe via \eqref{tragrad}\ that one has a canonical embedding of
$\ce_{[a_1\cdots a_k]}[k]=\La^k\ce_a[1]$ into $(\La^k\TT)^0/(\La^{k}\TT)^1=\gr_{0}(\La^k\TT)$. Since $\tau_-\in\TT[-1]$, $\si\mapsto \tau_-\wedge\si$ is thus seen to yield an isomorphism of $\ce_{a_1\cdots a_k}[k+1]$ with $\gr_{-1}(\La^{k+1}\TT)$ and analogously for the other components.

The tractor connection on $\La^{k+1}\TT$ is given by
\begin{gather}
     \na^{\La^{k+1}\TT}_c
    \begin{pmatrix}
            \rh_{a_1\cdots a_k} \\
            \ph_{a_0\cdots a_k}\; |\;\ \mu_{a_2\cdots a_k} \\
            \si_{a_1\cdots a_k}
    \end{pmatrix}
    \nonumber\\ \qquad{}=
\left(
\begin{matrix}
      {D}_c\rh_{a_1\cdots a_k} -P_{c}^{\; p}\ph_{pa_1\cdots a_k}-kP_{c[a_1}\mu_{a_2\cdots a_k]} \\
\left(
\begin{matrix}
      {D}_c\ph_{a_0\cdots a_k} +(k+1)\bg_{c[a_0}\rh_{a_1\cdots a_k]} \\
      +(k+1)P_{c[a_0}\si_{a_1\cdots a_k]}
\end{matrix}
\right)
      \; | \;\
\left(
\begin{matrix}
 {D}_c\mu_{a_2\cdots a_k}
      \\
      -P_{c}^{\; p}\si_{pa_2\cdots a_k}
      +\rh_{ca_2\cdots a_k}
\end{matrix}
\right)
      \\
      {D}_c\si_{a_1\cdots a_k}
      - \ph_{ca_1\cdots a_k}+k\bg_{c[a_1}\mu_{a_2\cdots a_k]}.
    \end{matrix}\label{tracon}
    \right).
\end{gather}
The f\/irst BGG-splitting operator $L_0^{\La^{k+1}\TT}:\ce_{[a_1\cdots a_k]}[k+1]\goesto\La^{k+1}\TT$\
is given by
\begin{gather}
 L_0^{\La^{k+1}\TT}(\si)\label{L0form} \\
 \qquad{}=
  \begin{pmatrix}\!
    \begin{pmatrix}
          -\frac{1}{n(k+1)}D^pD_p\si_{a_1\cdots a_k}
          +\frac{k}{n(k+1)}D^pD_{[a_1}\si_{|p|a_2\cdots a_k]}
          +\frac{k}{n(n-k+1)}D_{[a_1}D^p\si_{|p|a_2\cdots a_k]}
          \\
          +\frac{2k}{n}P^p_{\; [a_1}\si_{|p|a_2\cdots a_k]}
          -\frac{1}{n}J\si_{a_1\cdots a_k}
    \end{pmatrix}
    \\
    {D}_{[a_0}\si_{a_1\cdots a_k]} \; | \;
    -\frac{1}{n-k+1}\bg^{pq}{D}_p\si_{qa_2\cdots a_k}
    \\
    \si_{a_1\cdots a_k}\!
  \end{pmatrix}\!.\nonumber
\end{gather}

\subsection{Almost Einstein scales}\label{section-aEs}

The f\/irst splitting operator for the standard
tractor bundle is
\begin{gather}\label{splitStd}
  L_0^{\TT}: \ \ce[1] \goesto \Ga(\TT),\qquad
  \si \mapsto
  \begin{pmatrix}
    \frac{1}{n}(\Lapl-J)\si \\
    D\si \\
    \si
  \end{pmatrix}.
 \end{gather}
By \eqref{sttracon},
\begin{gather*}
  \na^{\TT}\circ L_0^{\TT} (\si)=
  \begin{pmatrix}
    \frac{1}{n}D_c(\Lapl\si-J\si)-P_c^{\; p}D_{p}\si  \\
    (D_aD_b\si+P_{ab}\si)+\frac{1}{n}(\Lapl\si-J\si)\bg_{ab}  \\
    0
  \end{pmatrix}.
\end{gather*}
Since $\frac{1}{n}(\Lapl\si-J\si)\bg_{ab}$ is minus the
trace-part of $(D_aD_b\si+P_{ab}\si)$ and $\HH_1(\TT)=\ce_{(ab)_0}$\
we have that the f\/irst BGG-operator of $\TT$ is
\begin{gather*}
  \Th_0^{\TT}: \ \ce[1] \goesto\ce_{(ab)_0}, \qquad
  \si \mapsto (D_aD_b\si+P_{ab}\si)_0.
\end{gather*}
By computing the change of the Schouten tensor $P$ with respect
to a conformal rescaling one obtains that with $U=\{x\in M:\si(x)\not=0\}$,
\begin{gather}\label{aEs-equ}
  (D_aD_b\si+P_{ab}\si)_0=0 \ \equiv \ \si^{-2}\bg\ \mr{is\ Einstein\ on}\ U.
\end{gather}
This says that $P(\si^{-2}\bg)$, or equivalently $\mr{Ric}(\si^{-2}\bg)$, is
a multiple of $\si^{-2}\bg$ on $U$.
$U$ always has to be an open dense subset of $M$,
and we call the set of solutions of \eqref{aEs-equ}
the space of \emph{almost Einstein scales}~\cite{gover-aes},
i.e.
\begin{gather}\label{def-aEs}
  \mb{aEs}([g])=\ker \Th_0^{\TT}\subset\ce[1].
\end{gather}
It turns out to be a dif\/ferential consequence of \eqref{aEs-equ}
that $\frac{1}{n}D_c(\Lapl\si-J\si)=P_c^{\; p}D_{p}\si$, and thus one has the
well known fact
\begin{proposition}\label{prop-aEs}
  $\na^{\TT}$-parallel sections of the standard tractor bundle are
  in {\rm 1:1}-correspondence with $\mb{aEs}([g])$.
\end{proposition}

\subsection{Conformal Killing forms}\label{section-cKf}

Via \eqref{tracon} and \eqref{L0form} one computes
that for $\si\in\ce_{[a_1\cdots a_k]}[k+1]$
the projection of $\na^{\La^{k+1}\TT}\circ L_0^{\La^{k+1}\TT}(\si)$
to the lowest slot $\ce_{c[a_1\cdots a_k]}[k+1]$ in $\Om^1(M,\La^{k+1}\TT)$
is given by
\begin{gather*}
  \si_{a_1\cdots a_k} \mapsto D_c \si_{a_1\cdots a_k}
  -D_{[a_0}\si_{a_1\cdots a_k]}
    -\frac{k}{n-k+1}\bg^{pq}{D}_p\si_{qa_2\cdots a_k}.
\end{gather*}
This is the projection of $\si_{a_1\cdots a_k}$ to
the highest weight part of $\ce_{c[a_1\cdots a_k]}[k+1]$ which
is formed by trace-free elements with trivial alternation, we write
\begin{gather*}
  \ce_{\{c[a_1\cdots a_k]\}_0}[k+1]:=\{\si_{a_1\cdots a_k}: 0=\si_{[ca_1\cdots a_k]}\ \mr{and}\ 0=\bg^{ca_1}\si_{ca_1\cdots a_k} \}.
\end{gather*}
One computes that in fact $\HH_1^{\La^{k+1}\TT}=\ce_{\{c[a_1\cdots a_k]\}_0}[k+1]$
and obtains the f\/irst BGG-operator
\begin{gather*}
  \Th_0^{\La^{k+1}\TT}: \ \ce_{[a_1\cdots a_k]}[k+1] \goesto
\ce_{\{c[a_1\cdots a_k]\}_0}[k+1],\qquad
  \si \mapsto D_{\{c}\si_{a_1\cdots a_k\}_0}.
\end{gather*}
Forms in the kernel of $\Th_0^{\La^{k+1}\TT}$
are thus the \emph{conformal  Killing}  $k$-\emph{forms}.

Unlike the case of $k=0$, it is not true for $k\geq 1$ that always $\na^{\La^{k+1}\TT}(L_0^{\La^{k+1}\TT}\si)=0$ for $\si\in\ker\Th_0^{\La^{k+1}\TT}\subset \ce_{[a_1\cdots a_k]}[k+1]$. However, given a section $s\in\Ga(\La^{k+1}\TT)$ with lowest slot
$\Pi(s)=\si\in\ce_{[a_1\cdots a_k]}$, one has by construction of $L_0^{\La^{k+1}\TT}$ that $s=L_0^{\La^{k+1}\TT}\si$ and that $\Th_0^{\La^{k+1}\TT}\si=\Pi\circ\na^{\La^{k+1}\TT}\circ L_0^{\La^{k+1}\TT}=0$; i.e., parallel sections of $\La^{k+1}\TT$\
do always project to special solutions of $\Th_0^{\La^{k+1}}\si=0$.
These solutions were termed \emph{normal}\ conformal Killing forms
by F.~Leitner~\cite{leitner-normal}.
Thus, by def\/inition, normal conformal Killing $k$-forms are
in 1:1-correspondence (via $L_0^{\La^{k+1}\TT}$ and $\Pi$)\
with $\na^{\La^{k+1}\TT}$-parallel sections of $\La^{k+1}\TT$.

Denote the components of the splitting $L_0^{\La^{k+1}\TT}\si$ given
in \eqref{L0form} by
\begin{gather}\label{splitdenom}
     \begin{pmatrix}
            \rh_{a_1\cdots a_k} \\
            \ph_{a_0\cdots a_k}\; |\;\ \mu_{a_2\cdots a_k} \\
            \si_{a_1\cdots a_k}
    \end{pmatrix}.
\end{gather}
A normal conformal
Killing form satisf\/ies $\na^{\La^{k+1}\TT}(L_0^{\La^{k+1}\TT}\si)$.
By~\eqref{tracon}, the resulting equation in lowest slot just says
that $\si$ is a conformal Killing form.
Additionally, we get the following equations for the components $\ph$, $\mu$, $\rh$:
\begin{gather}
{D}_c\rh_{a_1\cdots a_k} -P_{c}^{\; p}\ph_{pa_1\cdots a_k}-kP_{c[a_1}\mu_{a_2\cdots a_k]} =0,\nonumber \\
{D}_c\ph_{a_0\cdots a_k} +(k+1)\bg_{c[a_0}\rh_{a_1\cdots a_k]}
      +(k+1)P_{c[a_0}\si_{a_1\cdots a_k]}  =0, \label{normalcond}\\
 {D}_c\mu_{a_2\cdots a_k}
      -P_{c}^{\; p}\si_{pa_2\cdots a_k}
      +\rh_{ca_2\cdots a_k}  =0.\nonumber
\end{gather}

  \section{Generic rank two distributions\\  and  associated conformal structures}\label{rank2dist}

\subsection{The distributions}\label{sec-distributions}
Let $M$ be a 5-manifold. We are interested in generic  rank 2 distributions on $M$, i.e., rank 2 subbundles $\mathcal{D}\subset TM$  such that values of sections of $\mathcal{D}$ and Lie brackets of two such sections span  a rank 3 subbundle $[\mathcal{D},\mathcal{D}]$ and values of Lie brackets of at most three sections span the entire tangent bundle $TM$ (recall \eqref{distspan}). In other words, these are distributions of maximal growth vector $(2,3,5)$ in each point.

Def\/ining $T^{-1}M:=\mathcal{D}$, $T^{-2}M:=[\mathcal{D},\mathcal{D}]$ and $T^{-3}M=TM$,  the distribution gives rise to a~f\/iltration of the tangent bundle by subbundles compatible with the Lie bracket of vector f\/ields in the sense that for $\xi\in\Gamma(T^{i}M)$ and $\eta\in\Gamma(T^{j}M)$ we have $[\xi,\eta]\in\Gamma(T^{i+j}M)$.  Given such a f\/iltration, the Lie bracket of vector f\/ields induces a tensorial bracket $\mathcal{L}$ on the associated graded bundle $\mathrm{gr}(TM)=\bigoplus\mathrm{gr}_i(TM)$, where $\mathrm{gr}_i(TM)=T^{i}M/T^{i+1}M$. This bundle map $\mathcal{L}:\mathrm{gr}(TM)\times\mathrm{gr}(TM)\to\mathrm{gr}(TM)$ is called the \emph{Levi\ bracket}. It makes the bundle $\mathrm{gr}(TM)$ into a~bundle of nilpotent Lie algebras; the f\/iber $(\mathrm{gr}(TM)_x,\mathcal{L}_x)$ is the \emph{symbol algebra} at the point $x$. Note that a rank~2 distribution $\mathcal{D}\subset TM$ is generic if and only if the symbol algebra at each point is isomorphic to the graded Lie algebra $\g=\g_{-1}\oplus\g_{-2}\oplus\g_{-3}$,  where $\mathrm{dim}(\g_{-1})=2, $  $\mathrm{dim}(\g_{-2})=1, $ $\mathrm{dim}(\g_{-3})=2, $ and the only non-trivial components of the Lie bracket,  $\g_{-1}\times\g_{-2}\to\g_{-3}$ and $\Lambda^2\g_{-1}\to\g_{-2}$, def\/ine isomorphisms.

Generic rank 2 distributions in dimension~5 arise from ODEs of the form
\begin{gather}\label{ODE}
  z'=F(x,y,y',y'',z)
\end{gather}
with $\frac{\partial^2F}{\partial(y'')^2}\neq 0$ where $y$ and $z$ are functions of $x$, see e.g.~\cite{nurowski-metric} for that viewpoint. In his famous f\/ive-variables paper \cite{cartan-cinq} from 1910, \'{E}lie Cartan associated to these distributions a canonical Cartan connection, a result that shall be stated more precisely in~Section~\ref{fefferman}.

\subsection[Some algebra: $G_2$ in $\SO(3,4)$]{Some algebra: $\boldsymbol{G_2}$ in $\boldsymbol{\SO(3,4)}$}\label{somealgebra}

Let us recall one of the possible def\/initions of an exceptional Lie group of type $G_2$.
It is well known (e.g.~\cite{bryant-exceptional}) that the natural $\GL(7,\mathbb{R})$-action on the space $\Lambda^3{\rr^7}^*$ of 3-forms on $\rr^7$ has two open orbits and that the stabilizer of a 3-form in either of these open orbits is a 14-dimensional Lie group.  For one of these orbits it is   a compact real form of the complex exceptional Lie group~$G_2$, and for the other orbit it is a split real form. To distinguish between the two open orbits, consider the bilinear map $\rr^7\times\rr^7\to\Lambda^7{\rr^7}^*$, $(X,Y)\mapsto i_{X}\Phi\wedge i_{Y}\Phi\wedge\Phi$,
associated to a~3-form $\Phi$. This bilinear map  is non-degenerate if and only if $\Phi$ is contained in an open orbit. In that case, it determines  an invariant volume form $\mathrm{vol}$ on $\rr^7$ given by the root $^9\sqrt{D}\in\Lambda^7({\rr^7}^*)$ of its determinant $D\in(\Lambda^7{\rr^7}^*)^9$, see e.g.~\cite{hitchin-stable}. Hence
\begin{gather}\label{compatalg}
H(\Phi)(X,Y)\mathrm{vol}:=i_{X}\Phi\wedge i_{Y}\Phi\wedge\Phi
\end{gather} def\/ines  a $\mathbb{R}$-valued bilinear form $H(\Phi)$ on ${\rr^7}$ which is invariant under the action of the stabilizer of $\Phi$. It turns out that $H(\Phi)$  is positive def\/inite if the stabilizer is the compact real form, and it has signature $(3,4)$ if the stabilizer is the split real form of~$G_2$. In the sequel, let $G=G_2$ be the stabilizer of a 3-form $\Phi\in\Lambda^3{\rr^7}^*$ such that the associated bilinear form $H(\Phi)$ has signature~$(3,4)$.
The above discussion implies that this $G_2$ naturally includes into the special orthogonal group $\tilde{G}=\SO(3,4)$, an observation which will be crucial for what follows.

Let us be more explicit and realize $\SO(3,4)$ as $\SO(h)$, with
\begin{gather*}
  h=
  \begin{pmatrix}
      0 & 0 & 1 \\
  0 & M_{2,3} & 0 \\
  1 & 0 & 0
  \end{pmatrix},\qquad
  M_{2,3}=
  \begin{pmatrix}
   0 & 0 & 0 & 1 & 0 \\
   0 & 0 & 0 & 0 & 1 \\
   0 & 0 & -1 & 0 & 0  \\
   1 & 0 & 0 & 0 & 0  \\
   0 & 1 & 0 & 0 & 0
  \end{pmatrix}.
\end{gather*}
Via this bilinear form we identify ${\rr^7}^*\cong\rr^7$. Consider the standard basis $e_1,\dots, e_7$ on ${\rr^7}$, and def\/ine $G_2$ as the stabilizer of $\Phi\in\Lambda^3\rr^7,$
\begin{gather}
  \Phi:= -\tfrac{1}{\sqrt{3}} e_7\wedge e_2\wedge
e_3+\tfrac{1}{\sqrt{6}}e_5\wedge e_4\wedge
e_2+\tfrac{1}{\sqrt{6}}e_6\wedge e_4\wedge
e_3\nonumber\\
\phantom{\Phi:=}{} -\tfrac{1}{\sqrt{6}}e_7\wedge e_4\wedge e_1-\tfrac{1}{\sqrt{3}}e_1
\wedge e_5\wedge e_6.\label{Phialg}
\end{gather}
Then, via the identif\/ication $\Lambda^3{\rr^7}^*\cong\Lambda^3\rr^7$,
\begin{gather}\label{hphi}
H(\Phi)(X,Y)=\tfrac{1}{\sqrt{6}}h(X,Y),
\end{gather}
and this equation characterizes the $\SO(h)$-conjugacy class of $G_2$.
That is, $\Phi$\ has $\SO(h)$-stabilizer conjugated to $G_2$ if and only if $H(\Phi)$ is some non-zero multiple of~$h$.

The Lie algebra $\mathfrak{so}(h)$ has the matrix representation \eqref{solie}.
It contains the Lie algebra $\mathfrak{g}$ of $G_2$, which is formed by elements $M\in\mathfrak{gl}(7,\mathbb{R})$
 such that $ \Phi(M v,v',v'')+ \Phi(v,M v',v'')+ \Phi(v,v',M v'')=0$, as the subalgebra  consisting of matrices
\begin{gather}\label{g2matrices}
\begin{pmatrix}
    \mathrm{tr}(A)& \vline & Z &\vline &s&\vline & W& \vline& 0\\
    \hline
    X&\vline &A&\vline &\sqrt{2}\mathbb{J}Z^t&\vline&\frac{s}{\sqrt{2}}\mathbb{J}&\vline&-W^t\\
    \hline
    r&\vline&-\sqrt{2}X^t\mathbb{J}&\vline&0&\vline&-\sqrt{2}Z\mathbb{J}&\vline &s\\
    \hline
    Y&\vline&-\frac{r}{\sqrt{2}}\mathbb{J}&\vline &\sqrt{2}\mathbb{J}X&\vline &-A^t&\vline &-Z^t\\
    \hline
    0&\vline &-Y^t&\vline &r&\vline &-X^t&\vline &-\mathrm{tr}(A)
\end{pmatrix}
\end{gather}
with $A\in\mathfrak{gl}(2,\mathbb{R})$, $X,Y\in\rr^2$, $Z,W\in{\rr^2}^*$,
$r,s\in\mathbb{R}$ and $\mathbb{J}=\left(\begin{smallmatrix}0&-1\\
1&0\end{smallmatrix}\right)$.

For later use, let us note here that the complement of $\mathfrak{g}$ in $\mathfrak{so}(h)$ with respect to the Killing form is isomorphic to the seven dimensional standard representation of $G_2$. That means we have a  $G_2$-module decomposition
\begin{gather*}
\mathfrak{so}(h)=\mathfrak{g}\oplus{\rr^7}.
\end{gather*}
The sequence
  \begin{gather}\label{sequence-Phi}
    0\goesto\g\embed\so(h)\overset{i\Phi}{\goesto}\rr^7\goesto 0
  \end{gather}
  is $G_2$-equivariant and exact. Here
  \begin{gather*}
    i\Phi: \ \so(h)=\La^2\rr^7\goesto\rr^7
  \end{gather*}
  is the insertion of $\so(h)$ into $\Phi$.
The factor of $\Phi$ as given in \eqref{Phialg} was
chosen such that the insertion
  \begin{gather}\label{split2}
    i\Phi: \ \rr^7\goesto\La^2\rr^7=\so(h)
  \end{gather}
 splits sequence \eqref{sequence-Phi}.

Next we consider parabolic subgroups in $G_2$ and $\SO(h)$. Let $e_1\in{\rr^7}$ be the f\/irst basis vector in the standard representation. Then the isotropy group of the ray $\mathbb{R}_{+}e_1$ is a parabolic subgroup~$\tilde{P}$ in~$\SO(h)$, and the intersection $P=\tilde{P}\cap G_2$ is a parabolic subgroup in $G_2$.
To describe explicitly the corresponding parabolic subalgebra
$\mathfrak{p}\subset\mathfrak{g}$,
we introduce vector space decompositions of the Lie algebra.
We consider the block decomposition
 \[\begin{pmatrix}
   \mathfrak{g}_{0}&\mathfrak{g}_{1}&\mathfrak{g}_{2}&\mathfrak{g}_{3}&0\\
   \mathfrak{g}_{-1}&\mathfrak{g}_{0}&\mathfrak{g}_{1}&
    \mathfrak{g}_{2}&\mathfrak{g}_{3}\\
   \mathfrak{g}_{-2}&\mathfrak{g}_{-1}&0&\mathfrak{g}_{1}&\mathfrak{g}_{2}\\
   \mathfrak{g}_{-3}&\mathfrak{g}_{-2}&\mathfrak{g}_{-1}&\mathfrak{g}_{0}&
    \mathfrak{g}_{1}\\
   0&\mathfrak{g}_{-3}&\mathfrak{g}_{-2}&\mathfrak{g}_{-1}&\mathfrak{g}_{0}
 \end{pmatrix},
 \]
 of matrices \eqref{g2matrices}; this def\/ines a grading
 \[\mathfrak{g}=\mathfrak{g}_{-3}\oplus\mathfrak{g}_{-2}\oplus\mathfrak{g}_{-1}\oplus\mathfrak{g}_{0}\oplus\mathfrak{g}_{1}\oplus\mathfrak{g}_{2}\oplus\mathfrak{g}_{3}.\]
 Note that the subalgebra $\g_{-}=\g_{-3}\oplus\g_{-2}\oplus\g_{-1}$ coincides with the symbol algebra of a generic rank two distribution in dimension f\/ive as explained in Section~\ref{sec-distributions}.
 The grading  induces a~f\/iltration $\mathfrak{g}^3\subset\mathfrak{g}^2\subset\mathfrak{g}^1\subset\mathfrak{g}^0\subset\mathfrak{g}^{-1}\subset\mathfrak{g}^{-2}\subset\mathfrak{g}^{-3}$, which is preserved by the action of $P$ on $\mathfrak{g}$. The subalgebra $\mathfrak{p}=\mathfrak{g}^0$ is the Lie algebra of the parabolic $P$, and the subalgebra $\mathfrak{g}_0\cong\mathfrak{gl}(2,\mathbb{R})$ is the Lie algebra of the subgroup $G_0\subset P$ that even preserves the grading. The subgroup $G_0$ is  isomorphic to $\GL_+(2,\mathbb{R})=\{M\in \GL(2,\mathbb{R}):\mathrm{det}(M)>0\}$.

\subsection{The homogeneous model and associated Cartan geometries}\label{hommod}
Let us look at the Lie group quotient $G_2/P$ next.
The action of $G_2$ on the class $\mathrm{e}\tilde{P}\in \SO(h)/\tilde{P}$
induces a smooth map
\[G_2/P\to \SO(h)/\tilde{P}.\]
 Since both homogeneous spaces have the same dimension, the map is  an open embedding. Since $G_2/P$ is a quotient of a semisimple Lie group by a parabolic subgroup, it is compact, and the map is in fact a dif\/feomorphism.
 The group  $\SO(h)$ acts transitively on the space of null-rays in~$\rr^{7}$, which can be identif\/ied with the pseudo-sphere $Q_{2,3}\cong S^2\times S^3$. It turns out that the metric $h$
 on $\rr^7$ def\/ined in~\eqref{metrich} induces the conformal class of $(g_2,-g_3)$ on $Q_{2,3}$,
with $g_2$, $g_3$ being the round metrics on $S^2$ respectively $S^3$.
The pullback of that  conformal structure yields a~$G_2$-invariant conformal structure on~$G_2/P$.

Explicit descriptions of the canonical rank two distribution on $Q_{2,3}\cong S^2\times S^3$ can be found in~\cite{Sag06}. In an algebraic picture the distribution corresponds to the $P$-invariant subspace \mbox{$\mathfrak{g}^{-1}/\mathfrak{p}\subset\mathfrak{g}/\mathfrak{p}$}. Via the identif\/ication of $T(G/P)$ with $G\times_{P}\g/\p$, this invariant subspace gives rise to a~rank two distribution, which is generic in the sense of~Section~\ref{sec-distributions}.

More generally, suppose $(\mathcal{G},\omega)$ is any parabolic geometry
of type $(G_2,P)$. Recall from Section~\ref{section-general} that the Cartan connection $\omega$ def\/ines an isomorphism $TM\cong\mathcal{G}\times_{P}\mathfrak{g}/\mathfrak{p}$. Hence, for any such geometry, the subspace $\mathfrak{g}^{-1}/\mathfrak{p}$ gives rise to a rank two distribution. This distribution will  be generic if a regularity condition on the Cartan connection is assumed; we shall introduce this condition next:
Let
\[T^{-1}M\subset T^{-2}M\subset TM\] be the sequence of subbundles of constant ranks $2$, $3$ and $5$ coming from the $P$-invariant f\/iltration  $\g^{-1}/\p\subset\g^{-2}/\p\subset\g/\p$. Consider the associated graded bundle $\mathrm{gr}(TM).$ This bundle can be naturally identif\/ied with $\mathcal{G}\times_P\mathrm{gr}(\g/\p)\cong(\mathcal{G}/P_+)\times_{G_0}\g_{-}$.
Since the Lie bracket on the nilpotent Lie algebra $\g_{-}$ is
invariant under the $G_0$-representation,
it induces a bundle map $\{,\}:\mathrm{gr}(TM)\times\mathrm{gr}(TM)\to\mathrm{gr}(TM)$, the \emph{algebraic bracket}. A Cartan connection form $\omega$ is said to be \emph{regular} if the f\/iltration $T^{-1}M\subset T^{-2}M\subset TM$ is compatible with the Lie bracket of vector f\/ields and the algebraic bracket coincides with the Levi bracket of the f\/iltration. But this precisely means that the rank two subbundle $\mathcal{D}:=T^{-1}M$ is generic and $T^{-2}M=[\mathcal{D},\mathcal{D}]$ (compare with Section~\ref{sec-distributions} and the structure of $\g_{-}$).
Regularity can be expressed as a condition on the curvature of a Cartan connection.
Since $\mathfrak{g}$ has a $P$-invariant f\/iltration, we have a  notion of maps in $\Lambda^k(\g/\p)^*\otimes\g$ of homogeneous degree $\geq l$, and the set of these maps is $P$-invariant.
A Cartan connection form is  regular if and only if the curvature function is homogeneous of degree $\geq 1$; this means that $\kappa(u)(\g^{i},\g^{j})\subset\g^{i+j+1}$ for all $i$, $j$ and $u\in\mathcal{G}$. Note that if the curvature function takes values in $\Lambda^2(\g/\p)^*\otimes\p$, i.e.\ the geometry is \emph{torsion-free}, then it is regular.

Now we can state Cartan's classical result in modern language. In this paper we restrict our considerations to orientable distributions. Equivalently,  this means that the bundle $TM$ be orientable. Then we have the following:

\begin{theorem}[\cite{cartan-cinq}]\label{thm-rank2}
One can naturally associate a regular, normal parabolic geometry $(\mathcal{G},\omega)$
of type $(G_2,P)$ to an  orientable generic rank two distribution in dimension five,  and this establishes an equivalence of categories.
\end{theorem}
The above discussion explains  that a regular parabolic geometry of type $(G_2,P)$ determines an underlying generic rank two distribution $\mathcal{D}$, and (for our choice of $P$) the  distribution turns out to be orientable.
The converse is shown in two steps: First one constructs a regular parabolic geometry inducing the given distribution. Next one employs an inductive normalization procedure based on the proposition below, which we state explicitly here, since it will be needed it in Proposition \ref{normaluptohom}.
\begin{proposition}[\cite{cap-slovak-par}]\label{prop-normdeform}
Let $(\mathcal{G},\omega)$ be a regular parabolic geometry with curvature function $\kappa$, and suppose that $\partial^*\kappa$ is of homogeneous degree $\geq l$ for some $l\geq 1$.  Then, there is a normal Cartan connection $\omega_{N}\in\Omega^1(\mathcal{G},\mathfrak{g})$ such that  $(\omega_{N}-\omega)$ is of homogeneous degree $\geq l$.
\end{proposition}
In the proposition the dif\/ference $\omega_N-\omega$, which
is horizontal, is viewed as a function $\mathcal{G}\to(\g/\p)^*\otimes\g$,
and the homogeneity condition employs the canonical f\/iltration of $(\g/\p)^*\otimes\g$.

\subsection{A Fef\/ferman-type construction}\label{fefferman}

The relation in Section~\ref{hommod} between the homogeneous models $G_2/P$ and $\SO(h)/\tilde{P}$
 suggests a relation between Cartan geometries of type $(G_2,P)$ and $(\SO(h),\tilde{P})$, i.e., generic rank two distributions and conformal structures.
Indeed, it was P.~Nurowski who f\/irst observed in \cite{nurowski-metric} that any generic rank two distribution on a f\/ive manifold $M$ naturally determines a conformal class of metrics of signature $(2,3)$ on $M$. Starting from a system \eqref{ODE} of ODEs, he explicitly constructed a metric from the conformal class. A dif\/ferent construction of such a metric can be found in \cite{cap-sagerschnig-metric}.

In the present paper,  we shall discuss Nurowski's result as a special case of an extension functor of Cartan geometries, see \cite{cap-zadnik-chains,cap-constructions,doubrov-slovak-inclusions}. Let $i':\mathfrak{g}\hookrightarrow\so(h)$  denote the derivative of the inclusion  $i:G_2\hookrightarrow \SO(h)$. Given a Cartan geometry $(\mathcal{G}\to M,\omega)$ of type $(G_2,P)$, we can extend the structure group of the Cartan bundle, i.e., we  can form the associated bundle $\tilde{\mathcal{G}}=\mathcal{G}\times_{P}\tilde{P}$. Then this is a principal bundle over $M$ with structure group $\tilde{P}$.
We have a natural inclusion $j:\mathcal{G}\hookrightarrow\tilde{\mathcal{G}}$ mapping  an element $u\in\mathcal{G}$ to the class  $[(u,e)]$. Moreover, we  can uniquely extend the Cartan connection $\omega$ on $\mathcal{G}$ to a Cartan connection $\tilde{\omega}\in\Omega^1(\tilde{\mathcal{G}},\mathfrak{so}(h))$
such that $j^*\tilde{\omega}=i'\circ \omega$.
The construction indeed def\/ines a functor from Cartan geometries of type $(G_2,P)$ to Cartan geometries of type $(\tilde{G},\tilde{P})$.

We will later need the relation between the curvatures
of $\ti\om$ and $\om$, which is discussed
  in the next lemma. We use the inclusion of
  adjoint tractor bundles $\A M\embed\ti\A M$ via
  \begin{gather*}
    \A M=\G\times_P \g\embed \G\times_P \so(h)=\ti\G\times_{\ti P}\so(h)=\ti A M.
  \end{gather*}
  \begin{lemma}\label{lem-curvaturerelations}\mbox{}
    \begin{enumerate}\itemsep=0pt
    \item[$1.$]     The curvature form $\ti\Om\in\Om^2(\ti\G,\so(h))$ of $\ti\om$\
    pulls back to the curvature form $\Om\in\Om^2(\G,\g)$ of $\om$:
    \begin{gather*}
      j^*(\ti\Om)=\Om.
    \end{gather*}
  \item[$2.$] The factorizations $K\in\Om^2(M, \A M)$ of $\Om$ and $\ti K\in\Om^2(M,\ti \A M)$ of $\ti\Om$ agree:
  \begin{gather*}
    \ti K=K \in\Om^2(M,\A M).
  \end{gather*}
    \end{enumerate}
    In particular, $K\in\Om^2(M,\ti\A M)$
  has values in $\A M\subset\ti\A M$.
 \end{lemma}

  \begin{proof}
      Since the exterior derivative $d$ is natural, it commutes
  with pullbacks: $j^*d\ti\om=d(j^*\ti\om)=d\om$.
  Since also $j^*([\ti\om,\ti\om])=[\om,\om]$,
  we thus see that by def\/inition of curvature \eqref{def-curvature}
  we have $j^* {\ti\Om}={\Om}$.

  Now the inclusion $j:\G\goesto\ti\G$ is a reduction of
  structure group from $\ti P$ to $P$. Therefore,
  factorizing $\ti{\Om}\in\Om^2_{\mr{hor}}(\ti\G,\so(h))^{\ti P}$\
  to the curvature form $\ti K\in\Om^2(M,\ti\A M)$ is the same
  as pulling back $\ti{\Om}$ via $j$ and then factorizing.
  \end{proof}

By Theorem \ref{thm-rank2}, we can associate a canonical Cartan geometry $(\mathcal{G},\omega)$ of type $(G_2,P)$ to a~generic rank two distribution on a~f\/ive manifold $M$. As discussed in Section~\ref{section-conformal}, any  Cartan geometry $(\tilde{\mathcal{G}},\tilde{\omega})$ of type $(\tilde{G},\tilde{P})$ determines a conformal structure on the underlying manifold $M$. Thus the above Fef\/ferman construction shows that a generic rank two distribution $\mathcal{D}$ naturally determines a conformal class $[g]$ of metrics of signature $(2,3)$. However, a priori we do not know whether $\tilde{\omega}$ is the normal Cartan connection associated to that conformal structure (which will be important for applications); to see this requires a proof. We will give a proof based on the following result which is derived via  BGG-techniques:
 \begin{proposition}[\cite{cap-correspondence}]\label{E}
Suppose $\mathbb{E}\subset\mathrm{ker}(\partial^*)\subset\Lambda^2(\mathfrak{g}/\mathfrak{p})^*\otimes\mathfrak{g}$ is a $P$-submodule, and consider the $G_0$-module $\mathbb{E}_0:=\mathbb{E}\cap\mathrm{ker}(\square). $
Let $(\mathcal{G}\to M,\omega)$ be a regular, normal parabolic geometry such that the harmonic curvature takes values in $\mathbb{E}_0$. If either $\omega$ is torsion-free or for any $\phi,\psi\in\mathbb{E}$ we have $\partial^*(i_{\phi}\psi)\in\mathbb{E}$, where $i_{\phi}\psi$ is the alternation of the map  $(X_0,X_1,X_2)\mapsto\phi(\psi(X_0,X_1)+\mathfrak{p},X_2)$ for $X_i\in\mathfrak{g}/\mathfrak{p}$,  the curvature function $\kappa$ has takes in $\mathbb{E}$.
 \end{proposition}

Kostant's version of the Bott--Borel--Weil theorem \cite{kostant-61} provides an algorithmic description of the $G_0$-representation $\mathrm{ker}(\square)$. Doing the necessary calculations for the $|3|$-graded Lie algebra~$\mathfrak{g}$ of~$G_2$ from Section~\ref{somealgebra} (or consulting~\cite{silhan-real}, if you want the computer to calculate for you) leads~to:

\begin{lemma} \label{lem-torsionfree}
The harmonic component $\kappa_{H}$ of a regular, normal parabolic geometry of type $(G_2,P)$ takes values  in a  $G_0=\GL_+(2,\mathbb{R})$-submodule of $\mathfrak{g}_1\wedge\mathfrak{g}_3\otimes\mathfrak{g}_0$  isomorphic to $S^4(\mathbb{R}^2)^*$.  It follows, that the geometry is torsion-free.
\begin{proof}
An algorithmic calculation shows that the  $G_0=\GL_+(2,\mathbb{R})$-representation $\mathrm{ker}(\square)$ is an irreducible summand of $\mathfrak{g}_1\wedge\mathfrak{g}_3\otimes\mathfrak{g}_0$  isomorphic to $S^4(\mathbb{R}^2)^*$.
Since $\kappa_{H}$ takes values in $\mathrm{ker}(\square)$, it is, in particular, contained in $\Lambda^2\mathfrak{p}_+\otimes\mathfrak{p}$. Now we can apply  Proposition~\ref{E} to conclude that the entire curvature $\kappa$ takes values in that $P$-module, i.e., the geometry is  torsion-free.
\end{proof}
\end{lemma}

We now show that
normality of $\om$ implies normality of $\ti\om$:
\begin{proposition}\label{normality}
Let $(\mathcal{G}\to M,\omega)$ be a regular normal parabolic geometry of type $(G_2,P)$, and let $(\tilde{\mathcal{G}}\to M,\tilde{\omega})$ be the associated parabolic geometry of type $(\tilde{G},\tilde{P})$. Then $\tilde{\omega}\in\Omega^1(\tilde{\mathcal{G}},\so(h))$ is normal.
\end{proposition}

\begin{proof}
The Killing forms of $\g$ resp.\ $\so(h)$ provide natural identif\/ications $\mathfrak{p}_{+}\cong(\mathfrak{g}/\mathfrak{p})^*$ and $\tilde{\mathfrak{p}}_{+}\cong(\so(h)/\tilde{\mathfrak{p}})^*$.
This allows us to view the curvature function of the geometry $(\G,\om)$
as a function $\kappa:\mathcal{G}\to\Lambda^2\mathfrak{p}_{+}\otimes\mathfrak{g}$ and the curvature function of $(\ti\G,\ti\om)$ as a function
$\tilde{\kappa}:\tilde{\mathcal{G}}\to\Lambda^2\tilde{\mathfrak{p}}_{+}\otimes\so(h)$.
The inclusion $i':\mathfrak{g}\to\so(h) $ and the isomorphism $\mathfrak{g}/\mathfrak{p}\cong\so(h)/\tilde{\mathfrak{p}}$ induce a map $I:\Lambda^2\mathfrak{p}_{+}\otimes\mathfrak{g}\to\Lambda^2\tilde{\mathfrak{p}}_{+}\otimes\so(h)$.  Since all  maps involved in the construction are homomorphisms of $P$-modules, this is indeed a homomorphism of $P$-modules. In terms of  $I$, the curvature functions are related via
\begin{gather}\label{curvrel}
\tilde{\kappa}(j(u))=I\circ\kappa(u)
\end{gather}
for all $u\in\mathcal{G}$. By equivariance, this uniquely determines $\tilde{\kappa}$ on $\tilde{\mathcal{G}}$.

Let $\tilde{\partial}^*:\Lambda^2\tilde{\mathfrak{p}}_{+}\otimes\so(h)\to\tilde{\mathfrak{p}}_{+}\otimes\so(h)$ be the Kostant codif\/ferential describing the conformal normalization condition, i.e.,
 on decomposable elements $U\wedge V\otimes A\in \Lambda^2\tilde{\mathfrak{p}}_{+}\otimes\so(h)$ we have
 \begin{gather}\label{del*}
\tilde{\partial}^*(U\wedge V\otimes A)= U\otimes [V,A]-V\otimes [U,A].
\end{gather}
To prove that the geometry $(\tilde{\mathcal{G}},\tilde{\omega})$ is normal amounts to showing that  $\tilde{\kappa}$ takes values in $\mathrm{ker}(\tilde{\partial}^*)$. By $\eqref{curvrel}$, this is equivalent to the fact that $\kappa$ takes values in $\mathrm{ker}(\tilde{\partial}^*\circ I)$. Proposition \ref{E}\  allows to reduce this remaining problem to a purely algebraic one:
Since, by Lemma \ref{lem-torsionfree}, regular, normal parabolic geometries of type $(G_2,P)$ are torsion-free,  Proposition \ref{E}  implies that $\kappa$ takes values in $\mathrm{ker}(\tilde{\partial}^*\circ I)$ if and only if this is true for the harmonic curvature component $\kappa_{H}$. By equivariance,  the map $\tilde{\partial}^*\circ I$ either vanishes on  $G_0$-irreducible components, or it is an isomorphism.  We have observed that $\mathrm{ker}(\square)$ is an irreducible $G_0$-representation isomorphic to $S^4(\mathbb{R}^2)^*$. Hence if we can show that the image of the map $\tilde{\partial}^*\circ I:\Lambda^2\mathfrak{p}_{+}\otimes\mathfrak{g}_0\to \tilde{\mathfrak{p}}_{+}\otimes\so(h)$ does  not contain an irreducible summand isomorphic to $S^4(\mathbb{R}^2)^*$, then $\kappa_{H}$ has to be contained in that kernel. But looking at  formula \eqref{del*}, we see that $\tilde{\partial}^*\circ I(\mathfrak{p}_{+}\otimes\mathfrak{g}_0)$  is actually contained in $\tilde{\mathfrak{p}}_{+}\otimes\tilde{\mathfrak{p}}_{+},$ which is easily seen to contain no summand isomorphic to $S^4(\mathbb{R}^2)^*$.
\end{proof}

\subsection{The parallel tractor three-form\\ and the underlying conformal Killing 2-form}

Let $\mathcal{T}$ be the standard tractor bundle for a conformal structure $[g]$ associated to a generic $2$-distribution $\mathcal{D}$ on a $5$-manifold $M$. Then $\mathcal{T}$ is easily seen to carry additional structure:

\begin{proposition}\label{prop2form}\mbox{}
\begin{enumerate}\itemsep=0pt
\item[$1.$]
The standard tractor bundle  $\mathcal{T}$ for a conformal structure $[g]$ associated to a generic $2$-distribution carries a parallel tractor $3$-form $\mathbf{\Phi}\in\Gamma(\Lambda^3\mathcal{T}^*)=\Gamma(\Lambda^3\mathcal{T})$.
\item[$2.$]
 The tractor $3$-form $\mathbf{\Phi}$ determines an underlying normal conformal Killing $2$-form $\phi\in\Gamma(\La^2T^*M\t\ce[3])=\ce_{[ab]}[3]$, which is locally decomposable.
\end{enumerate}
\end{proposition}

\begin{proof}
1.\ By construction, the conformal Cartan bundle is the associated bundle $\tilde{\mathcal{G}}=\mathcal{G}\times_{P}\tilde{P}$, where $\mathcal{G}$ is the Cartan bundle for  the distribution. Hence the tractor bundle can be viewed as $\mathcal{T}=\mathcal{G}\times_{P}{\rr^7}.$ It follows that the  $P$-equivariant function $f_{\Phi}:\mathcal{G}\to \Lambda^3{\rr^7}$ mapping constantly onto the three-form $\Phi$ stabilized by $G_2$ induces a  section $\mathbf{\Phi}\in\Gamma(\Lambda^3\mathcal{T})$. Proposition~\ref{normality} implies that the normal tractor connection $\nabla^{\Lambda^3\mathcal{T}}$ is induced from the normal Cartan connection $\omega\in\Omega^1(\mathcal{G},\mathfrak{g})$. Hence, according to \eqref{deftracon}, $\nabla^{\Lambda^3\mathcal{T}}_{\xi} \Phi$ corresponds to the function $u\mapsto (\xi'_u\cdot f_{\Phi})+\omega_u(\xi'(u))(f_{\Phi}(u))$, where $u\in\mathcal{G}$ and  $\xi'\in\mathfrak{X}(\mathcal{G})$ is a $P$-invariant lift of a the vector f\/ield $\xi\in\mathfrak{X}(M)$. Since $f_{\Phi}$ is constant and $\omega$ takes values in the isotropy algebra $\mathfrak{g}$ of $\Phi$, this means that  $\nabla^{\Lambda^3\mathcal{T}}\mathbf{\Phi}=0$, i.e.,  $\mathbf{\Phi}$ is a parallel tractor three-form.

2.\   Recall from Section~\ref{section-cKf} that we have a natural bundle projection  $\Pi_0:\Lambda^{3}\mathcal{T}\to\Lambda^2 T^*M\otimes\ce[3]$ and that  parallel sections of $\VV:=\Lambda^{3}\mathcal{T}$ project to normal conformal Killing $2$-forms.

Let $V=\La^3 \rr^7$  be endowed with
the canonical $\ti P$-invariant f\/iltration discussed in Section~\ref{conftra}.
We see from our explicit formula~\eqref{Phialg} for $\Phi$
that, up to a factor, a representative in $V/V^0=\gr_{-1}(V)$ is given by
 $e_7\wedge e_2\wedge e_3$. Around every point $x\in M$ we can choose a local section $\sigma:U\to\mathcal{G}$. On $\sigma(U)\subset\G\subset\ti\G$ we def\/ine $3$ constant functions, mapping  to~$e_7$, $e_2$ and $e_3$; these correspond to sections~$s_7$, $s_2$ and $s_3$ of the standard tractor bundle $\TT$.
$s_7$ is simply $\tau_-$;
to be precise, we use that $\si:U\goesto\G$ gives in particular a trivialization
of the conformal weight bundles, and we can regard $\tau_-$ as an (unweighted) section of $\TT$.
The tractors $s_2$ and $s_3$ lie in~$\TT^0$ and therefore
project to elements $\ph_2$ and $\ph_3$ in $\TT^0/\TT^1=\gr_0(\TT)=\ce_{a}$,
where we again use the trivialization of the conformal weight bundles.
Thus, $\tau_-\wedge\ph_2\wedge\ph_3$ is a representative of $\ph=\Pi_0(\mb{\Phi})\in \VV/\VV^0=\gr_{-1}(\VV)$, and the identif\/ication \eqref{identification} of $\VV/\VV^0=\gr_{-1}(\VV)$
with $\ce_{[ab]}$ tells us that $\ph=\ph_2\wedge\ph_3\in\ce_{[ab]}$.
\end{proof}

\begin{remark}
The parabolic subgroup $P$ of $G_2$ preserves a f\/iner f\/iltration of the standard representation, which yields the following ref\/inement of the 'conformal' f\/iltration
of standard tractor bundle:
\[\mathcal{T}^{2'}\subset\mathcal{T}^{1'}\subset\mathcal{T}^{0'}\subset\mathcal{T}^{-1'}\subset\mathcal{T}^{-2'}.\]
The isotropic line bundle $\mathcal{T}^{2'}$ corresponds to the subspace generated by $e_1\in{\rr^7}$. The bundle $\mathcal{T}^{-1'}$ is the orthogonal complement to $\mathcal{T}^{2'}$ with respect to the tractor metric.  The explicit form of the three-form $\Phi$ (see \eqref{Phialg}) shows how to characterize the additional f\/iltration components in terms of $\mathbf{\Phi}$. Let $\tau_+\in\Ga(\TT\t\ce[1])$ be the canonical insertion of $\ce[-1]$ into $\mathcal{T}^{2'}\subset \TT$ which
was already def\/ined in Section~\ref{conform-parabolic}.
 Then the subbundle $\mathcal{T}^{1'}$ can be described as the set of all tractors $s\in\Gamma(\mathcal{T})$ such that
$i_{s}i_{\tau_+}\mathbf{\Phi}=0.$
The subbundle $\mathcal{T}^{0'}$ is the bundle orthogonal  to $\mathcal{T}^{-1'}$.
\end{remark}

\begin{remark}
The distribution $\mathcal{D}$ can be recovered from the conformal class $[g]$ associated to the distribution  and the conformal Killing 2-form $\phi$. The kernel of the 2-form is the rank three distribution $[\mathcal{D},\mathcal{D}]$. The rank two distribution can be recovered as the kernel of the restriction of a metric $g$  to the rank three distribution, in other words, the set of isotropic elements in the kernel of $\phi$.
\end{remark}

\section{Holonomy reduction and characterization}\label{section-characterization}

The aim of the following section is to characterize conformal structures arising from generic rank~2 distributions in dimension f\/ive in terms of normal conformal Killing 2-forms satisfying certain additional equations. See Theorem A for a precise statement of the result.   We proceed as follows. First, we prove that a conformal manifold of signature $(2,3)$ whose conformal holonomy is contained in $G_2$ is obtained from a generic rank two distribution via a Fef\/ferman construction.  Then we aim for a  characterization of the conformal structures in terms of underlying conformal data;  we derive conditions to distinguish those  normal conformal Killing 2-forms  coming from parallel tractor 3-forms def\/ining holonomy reductions to $G_2$. This is done analogously to~\cite{cap-gover-holonomy-cr}, where the authors arrive at a version of Sparling's characterization \cite{graham-sparling}\ of Fef\/ferman spaces in terms of a conformal Killing f\/ield.

\begin{remark}
  Recently T.~Leistner and P.~Nurowski showed on examples that in some cases the conformal structures constructed in this way have explicit ambient metrics with holonomy group~$G_2$, see~\cite{nurowski-explicit} and~\cite{leistner-nurowski-g2ambient}.
\end{remark}

\subsection{Conformal holonomy}

Let $(M,[g])$ be a conformal structure of signature $(2,3)$
encoded in a Cartan geometry $(\ti\G,\ti\om)$ as described
in Section~\ref{section-conformal}.
The standard tractor bundle $\TT$ of $[g]$ is
endowed with the tractor connection $\na^{\TT}$, and we def\/ine the conformal holonomy of $[g]$ as
\begin{gather*}
  \mr{Hol}([g]):=\mr{Hol}(\na^{\TT}).
\end{gather*}
See also \cite{armstrong-conformal}.
Now $\TT$ comes about as associated bundle to $\ti\G':=\ti\G\times_{\ti P} \SO(h)$
and $\na^{\TT}$ is induced  from the principal
connection form $\ti\om'\in\Om^1(\ti\G',\so(h))$. Thus,
we have that $\mr{Hol}(\na^{\TT})=\mr{Hol}(\ti\om')$.
By construction, the pullback of $\ti\om'$ to $\ti\G\subset\ti\G'$
is  the Cartan connection form~$\ti\om$.

In the Fef\/ferman-type construction of Section~\ref{fefferman}
we started with a parabolic geometry $(\G,\om)$ of type
$(G_2,P)$ encoding a generic rank $2$ distribution $\D$, and we associated
to this the parabolic geometry $(\ti\G,\ti\om)$
of type $(\SO(h),\ti P)$ by equivariantly extending $\om$ to $\ti\om$.
If we add the extended bundles $\G'=\G\times_P G_2$ and
$\ti\G'=\ti\G\times_{\ti P} \SO(h)=\G\times_P\SO(h)$ to the
picture, we obtain the commuting diagram of inclusions
\begin{gather*}
  \xymatrix{
    (\G',\om') \ar@{^{(}->}[r] & (\ti\G',\ti\om') \\
    (\G,\om) \ar@{^{(}->}[u] \ar@{^{(}->}[r] & (\ti \G,\ti\om) \ar@{^{(}->}[u]
  }
\end{gather*}
In particular, this yields a holonomy reduction of $(\ti\G',\ti\om')$
to $(\G',\om')$,
and thus $\mr{Hol}(\ti\om')=\mr{Hol}(\om')\subset G_2$. By Proposition~\ref{normality}, $\tilde{\omega}$ is normal. Hence, $\mr{Hol}(\ti\om')$ is indeed the conformal holonomy $\mr{Hol}([g]_{\D})$,  which is thus seen to be contained in~$G_2$.

We are now going to show the converse: if for a conformal
structure~$(M,[g])$ of signature~$(2,3)$ the conformal holonomy $\mr{Hol}([g])$ is contained in $G_2$,
then there is  a canonical generic rank $2$-distribution~$\D$
on $M$ such that $[g]=[g]_{\D}$.

Let $\pi:\ti\G'\goesto M$ be the surjective submersion of the
$\SO(h)$-principal bundle $\ti\G'$. The next proposition covers the holonomy
reduction of a conformal Cartan geometry to $G_2$. A similar result has been obtained in \cite{jesse-thesis}.

\begin{proposition}
\label{reduction}
  Let $(\ti\G,\ti\om)$ be such that $(\ti\G',\ti\om')$\
  has holonomy in $G_2$ and let $\HH\subset\ti\G'$ be a reduction
  of $(\ti\G',\ti\om')$ to $G_2$. Then:
  \begin{enumerate}\itemsep=0pt
  \item[$1.$] $\HH\subset\ti\G'$ and $\ti\G\subset\ti\G'$ intersect
    transversally. We denote the resulting submanifold by $\G:=\HH\cap\ti\G$.
  \item[$2.$] For every $u\in\G$, $T_u\pi(T_u\G)=T_{\pi(u)}M$.
  \item[$3.$] $\G$ is a $P$-principal bundle over $M$.
  \item[$4.$] Let $\om$ be the pullback of $\ti\om\in\Om^1(\G,\so(h))$\
    to $\G$. Then $\om\in\Om^1(\G,\g)$ is a Cartan connection form.
  \end{enumerate}
\end{proposition}

  \begin{proof} 1.\   We have that
$T_u\HH+T_u \ti\G\supset u\cdot \g+u\cdot\ti\p=u\cdot \so(h)=\ker(T_u\pi)$.
    Since $T_u\pi:T_u\ti\G\goesto T_{\pi(u)}M$ is surjective, we
    have that $\dim(T_u\HH+T_u \ti\G)=\dim (\so(h))+\dim(\g/\p)=\dim(T_u\ti\G')$.

2.\ Take $u\in\G=\HH\cap\ti\G$ and $\xi\in T_{\pi(u)}M$. Since
    the restrictions of $\pi$ to $\HH$ and $\ti\G$ are surjective
submersions, there exist $\xi_1\in T_u\HH$ and $\xi_2\in T_u\ti\G$\
such that $\xi=T_u\pi\xi_1=T_u\pi\xi_2$. Then
\begin{gather*}
  \xi_1-\xi_2\in\ker T_u\pi=u\cdot\so(h)=u\cdot \g+u\cdot\ti\p.
\end{gather*}
Thus, there exist $\eta_1\in u\cdot \g$ and $\eta_2\in u\cdot \ti\p$\
such that $\xi_1-\xi_2=\eta_1+\eta_2$. Let
\begin{gather*}
  \xi'=\xi_1-\eta_1=\xi_2+\eta_2\in T_u \G.
\end{gather*}
Then indeed $T_u\pi\xi'=\xi$.

3.\ Assume f\/irst that for an $x\in M$ there is a $u\in\HH_x\cap\ti\G_x=\G_x$.
  Then evidently $\G_x=u\cdot(G_2\cap\ti P)=u\cdot P$.
  It therefore remains to show that $\HH_x\cap\ti\G_x$ is always non-empty:
  let $u\in \HH_x$. Then there is a $g\in\SO(h)$ such that
  $u\cdot g\in\ti\G$. Since $G_2/P=\SO(h)/\ti P$ (see Section~\ref{somealgebra}),
  there is a $p\in\ti\P$ such that $gp=g'\in G_2$, and then $u\cdot g'\in\HH$\
  since $\HH$ is a $G_2$-subbundle and   $u\cdot g'=(u\cdot g)\cdot p\in\SO(h)$;
  i.e., $u\cdot g'\in \G$.

4.\ We now consider $\G$ as a reduction of the $\ti P$-principal bundle $\ti\G$ to $P$ and denote by $\om$ the pullback of $\ti\om\in\Om^1(\ti\G,\so(h))$. By construction, $\HH\subset\G$ was  obtained by holonomy reduction
  of $(\ti\G',\ti\om')$ to $G_2$. In particular, $\ti\om'_{|T\HH}$ has
  values in $\g$, and thus $\om\in \Om^1(\G,\g)$.
  $P$-equivariance and reproduction of $\p$-fundamental vector f\/ields
  is clear since $\G$ is just a $P$-principal subbundle of $\ti\G$\
  and $\ti\om$ is a Cartan connection form satisfying (C.\ref{cartan1})--(C.\ref{cartan2}) by assumption. We thus need to check that also (C.\ref{cartan3})
  holds for $\om$. I.e., for every $u\in\G$ we need that $\om_u:T_u\G\goesto \g$ is an isomorphism.
  We have seen that $T_u\pi(T_u\G)=T_{\pi(u)} M$.
  Since $u\cdot\ti\p=\ker(T_u\pi)\subset T_u\ti\G$ we see that
  $T_u\G\subset T_u \ti\G$ must span at least $\dim(\g/\p)$-complementary
  dimensions and thus already $T_u\G+u\cdot\p=T_u\ti\G$. But then
  $\om_u(T_u\G/u\cdot\p)=\ti\om_u(T_u\ti\g/u\cdot\ti\p)=\ti\g/\ti\p=\g/\p$.
  This, together with $\om_u(u\cdot\p)=\p$ by reproduction of fundamental vector f\/ields, gives that indeed $\om_u(T_u\G)=\g$.
  \end{proof}

Now suppose  $(\ti\G,\ti\om)$ is a normal parabolic geometry of type $(\ti G,\ti P)$ associated to $[g]$ with $\mr{Hol}([g])\subset G_2$,  and let $(\G,\om)$ be the parabolic geometry of type $(G_2,P)$ obtained via reduction as explained in Proposition \ref{reduction}. Since every normal conformal Cartan connection is torsion-free and $\mathfrak{p}\subset\tilde{\p}$,  $\omega$ is torsion-free. This evidently implies that $\omega$ is regular, and therefore it determines an underlying generic $2$-distribution. We will show in Theorem \ref{holred} that the canonical conformal structure associated to this distribution is $[g]$;
this will employ the following

\begin{lemma}\label{normaluptohom}
Let $(\ti\G,\ti\om)$ be a normal parabolic geometry of type $(\ti G,\ti P)$ with $\mr{Hol}(\ti\om')\subset G_2$,  and let~$(\G,\om)$ be the parabolic geometry of type $(G_2,P)$ obtained via reduction. Then  there is a~normal Cartan connection $\omega_{N}\in\Omega^1(\G,\g)$ such that the difference $(\omega_{N}-\omega)$ is of homogeneous  degree~$\geq 3$.
\end{lemma}

\begin{proof}
 Let $\kappa:\mathcal{G}\to\Lambda^2\mathfrak{p}_{+}\otimes\mathfrak{g}$ be the  curvature function of  $\omega$. We will verify  that $\kappa$ is of homogeneous degree $\geq 3$. Since the Kostant codif\/ferential preserves homogeneities, then also $\partial^*\kappa$ is of homogeneous degree $\geq 3$, and by Proposition \ref{prop-normdeform} this shows that there is a normal Cartan connection $\omega_N$ that dif\/fers from $\omega$ at most in  homogeneous degree $\geq 3$.

Torsion-freeness of $\omega$ means
that for any $u\in\mathcal{G}$,  $\kappa(u)$ is contained in $\Lambda^2\mathfrak{p}_{+}\otimes\mathfrak{p}$. Hence the only component of $\kappa(u)$ of homogeneous degree $<3$  that remains to be investigated is
\[\kappa_2(u)\in\mathfrak{g}_{1}\otimes\mathfrak{g}_{1}\otimes\mathfrak{g}_{0}.\] We show that this component vanishes as well.

Let $\tilde{\partial}^*:\Lambda^2\tilde{\mathfrak{p}}_{+}\otimes\so(h)\to\tilde{\mathfrak{p}}_{+}\otimes\so(h) $ be the conformal Kostant codif\/ferential. Choose linearly independent elements $X_1,X_2\in\mathfrak{g}_{-1}$, $X_3\in\mathfrak{g}_{-2}$ and $X_4, X_5 \in\mathfrak{g}_{-3};$ then these elements give a  basis of $\so(h)/\tilde{\mathfrak{p}}$. Consider the dual basis  $\tilde{Z}_1,\dots,\tilde{Z}_5\in\tilde{\mathfrak{p}}_{+}$ with respect to the Killing form on $\so(h)$. By construction, $\kappa(u)(X_i,X_j)=\tilde{\kappa}(u)(X_i,X_j)$ for all  $u\in\mathcal{G}.$ Thus normality of the conformal Cartan connection implies
\begin{gather}\label{gleichg1}
\sum_{i<j}([\tilde{Z}_i,\kappa(u)(X_i,X_j)]\otimes\tilde{Z}_j-[\tilde{Z}_j,\kappa(u)(X_i,X_j)]\otimes\tilde{Z}_i)=0.
\end{gather}

Recall that we have a $\mathfrak{g}$-module decomposition $\so(h)=\mathfrak{g}\oplus\mathbb{V}$. The projection $\pi_{\mathfrak{g}}:\so(h)\to\mathfrak{g}$ maps $\tilde{Z_j}$ to an element $Z_j\in\mathfrak{p}_{+}$ dual to $X_j\in\mathfrak{g}_{-}$ with respect to a multiple of the Killing form on~$\g$. Equivariance of the projection implies $\pi_{\mathfrak{g}}([\tilde{Z}_i,\kappa(u)(X_i,X_j)])=[Z_i,\kappa(u)(X_i,X_j)]$. It follows that $\eqref{gleichg1}$ also holds with $\tilde{Z}_i$'s replaced by $Z_i$'s.
The only part of that sum contained in lowest homogeneity, i.e. in $\mathfrak{g}_1\otimes\mathfrak{g}_1$, is
\begin{gather*}
[Z_1,\kappa_2(u)(X_1,X_2)]\otimes Z_2-[Z_2,\kappa_2(u)(X_1,X_2)]\otimes Z_1,
\end{gather*}
and so this expression vanishes as well.
Since the representation of $\mathfrak{g}_0$ on~$\mathfrak{g}_1$ given by the Lie bracket is faithful, this indeed implies that $\kappa_2(u)(X_1,X_2)=0$.
\end{proof}

Having Proposition \ref{reduction} and Lemma \ref{normaluptohom}, we can now show:

\begin{theorem}\label{holred}
Let $(M,[g])$ be a conformal structure of signature $(2,3)$ with conformal holonomy $\mr{Hol}([g])\subset G_2$.
Then $[g]$ is canonically associated to a generic rank two distribution $\mathcal{D}$ via a~Fefferman-type construction.
\end{theorem}

\begin{proof}
Let $(\ti\G,\ti\om)$ be the normal parabolic geometry of type $(\SO(h),\ti P)$ associated to the conformal structure $[g]$. Let $(\mathcal{G},\omega)$ be the parabolic geometry of type $(G_2,P)$ constructed in Proposition~\ref{reduction}.
Then we know that $\om$ is regular, and by Lemma~\ref{normaluptohom} there is a normal Cartan connection $\omega_{N}\in\Omega^1(\G,\g)$ that dif\/fers from $\omega$ at most in homogeneous degree $\geq 3$.

 Recall that the Cartan connection $\omega$ determines an isomorphism $\mathcal{G}\times_{P}\mathfrak{g}/\mathfrak{p}\cong TM$. Regularity of $\omega$ implies that the image of $\mathcal{G}\times_P\mathfrak{g}^{-1}/\mathfrak{p}$ under this isomorphism is a generic rank two distribution $\mathcal{D}$. Furthermore, we have a $P$-invariant conformal class of bilinear forms of signature $(2,3)$ on $\g/\p$, and the  conformal structure induced via the above isomorphism on $M$ is just $[g]$. On the other hand, the Fef\/ferman construction associates a conformal structure  $[g]_{\mathcal{D}}$ to the distribution~$\mathcal{D}$. This is the conformal structure induced via the isomorphism $\mathcal{G}\times_{P}\mathfrak{g}/\mathfrak{p}\cong TM$ def\/ined by the normal Cartan connection $\omega_{N}\in\Omega^1(\G,\g)$ associated to the distribution~$\mathcal{D}$. Since $\omega_{N}-\omega$ is of homogeneous degree $\geq 3$, the dif\/ference $(\omega-\omega_{N})$ takes values in $\mathfrak{p}$. But this implies that~$\omega$ and~$\omega_{N}$ induce the same isomorphism $TM\cong\mathcal{G}\times_{P}\mathfrak{g}/\mathfrak{p}$  and hence the same conformal structure on~$M$; i.e., the conformal
structure $[g]$ is the one induced by the distribution $\D$: $[g]=[g]_{\mathcal{D}}$.
\end{proof}

\subsection{Characterization via the tractor 3-form}\label{Phired}

We have seen that conformal structures associated to generic rank two distributions in dimension f\/ive precisely correspond to reductions in conformal holonomy from $\SO(h)$ to $G_2$. As a next step towards the desired characterization result, we explain that such a holonomy reduction can be encoded in terms of a parallel tractor $3$-form satisfying a certain compatibility condition with the tractor metric.

Let $\ti\G'$ be  the
extended $\SO(h)$-principal bundle over $M$, and let
 $\ti\om'\in\Om^1(\ti\G',\so(h))$ be the extension of the Cartan connection to a principal connection form.
We consider the  holonomy group $\Hol_u=\Hol_u(\ti\om')$ for an arbitrary point $u\in\ti\G'$.
Then, for $g\in \SO(h)$, one has $\Hol_{u\cdot g}=g\Hol_ug^{-1}$, and
$\Hol([g])$  is well def\/ined up to conjugation in $\SO(h)$.

\def\H{\mathcal{H}}

Let $\H_u\embed\ti\G'$, $u\in\ti\G'$, be the reduction of the $\SO(h)$-bundle
 $\ti\G'$ to $\Hol_u$. If $\mb{\Psi}\in\La^3\TT$ is parallel, it corresponds
to a $\SO(h)$-equivariant function $f:\ti\G'\goesto \La^3\rr^{7}$\
which is a constant $\Psi_u\in\La^3\rr^{7}$ on $\H_u$.
Hence $\Hol_u\cdot \Psi_u=\Psi_u$, or $\Hol_u\subset \SO(h)_{\Psi_u}$.
If $u'$ is another point in $\ti\G'$ one has $\Hol_{u'}=g\Hol_ug^{-1}$ for
some $g\in\SO(h)$ and $\Psi_{u'}=g\cdot \Psi_u$. Thus
$f(\ti\G')=\SO(h)\cdot\Psi_u$. We say that $\SO(h)_{\Psi_u}$\
is the \emph{orbit\ type}\
of the parallel tractor $\mb{\Psi}$. To be precise, the orbit type is def\/ined up to conjugation in $\SO(h)$.

Recall that the group $G_2\subset \SO(h)$ has been realized as the stabilizer
$\SO(h)_{\Phi}$ for $\Phi\in\La^3\rr^{7}$ given by \eqref{Phialg}.
Compatibility condition \eqref{hphi} singles out the $\SO(h)$-orbit
of $\Phi\in\La^3\rr^{7}$.
Hence $\Hol([g])$ reduces to $G_2$ if and only if there is a
$\na^{\La^3\TT}$-parallel $\mb{\Phi}\in\La^{3}\TT$ satisfying
the global version of \eqref{hphi}, i.e.,
\begin{gather}\label{compat}
  H(\mb{\Phi})=\la \mb{h}\ \mr{for\ a}\ \la\in\rr\backslash\{0\},
\end{gather}
where, for $s_1,s_2\in\Ga(\TT)$,
\begin{gather}\label{defH}
  H(\mb{\Phi})(s_1,s_2)\mb{vol}=i_{s_1}\mb{\Phi}\wedge i_{s_2}\mb{\Phi}\wedge\mb{\Phi},
\end{gather}
and $\mb{vol}\in \La^{7}\TT^*$ is the tractor volume form.

\subsection{Characterization in terms of the underlying conformal Killing 2-form}\label{2form}

We want to express compatibility condition~\eqref{compat}
of the $\na^{\La^3\TT}$-parallel tractor $\mb{\Phi}\in\La^3\TT$
in terms of the underlying normal conformal Killing 2-form $\phi=\Pi_0(\mb{\Phi})
\in \Ga(\La^2T^*M\t\ce[3])$.

According to \eqref{splitdenom}\ we have
\begin{gather*}
  \mb{\Phi}=
    \begin{pmatrix}
            \rh_{a_1 a_2} \\
            \ph_{a_0 a_1 a_2}\; |\;\ \mu_{a_2} \\
            \phi_{a_1 a_2}
    \end{pmatrix}
    =
  \begin{pmatrix}
    \begin{pmatrix}
          -\frac{1}{15}D^pD_p\phi_{a_1 a_2}
          +\frac{2}{15}D^pD_{[a_1}\phi_{|p|a_2]}
          +\frac{1}{10}D_{[a_1}D^p\phi_{|p|a_2]}
          \\
          +\frac{4}{5}P^p_{\; [a_1}\phi_{|p|a_2]}
          -\frac{1}{5}J\phi_{a_1a_2}
    \end{pmatrix}
    \\
    {D}_{[a_0}\phi_{a_1a_2]} \; | \;
    -\frac{1}{4}\bg^{pq}{D}_p\phi_{qa_2}
    \\
    \phi_{a_1 a_2}
  \end{pmatrix}.
\end{gather*}
We will also write this splitting  as
\begin{gather}\label{splitPhi}
    \Phi=\tau_-\wedge\phi+\ph+\tau_+\wedge\tau_-\wedge\mu+\tau_+\wedge\rh = \begin{pmatrix}
            \rh \\
            \ph  \; |\;\ \mu \\
            \phi
    \end{pmatrix}.
   \end{gather}
By \eqref{normalcond}\ the  conditions for
a conformal Killing $2$-form $\phi\in\Om^2(M)\t\ce[3]$\
to be normal are
\begin{gather*}
{D}_c\rh_{a_1 a_2} -P_{c}^{\; p}\ph_{pa_1 a_2}-2P_{c[a_1}\mu_{a_2]}=0, \\
{D}_c\ph_{a_0a_1a_2} +3\bg_{c[a_0}\rh_{a_1 a_2]}
      +3P_{c[a_0}\phi_{a_1 a_2]}  =0, \\
 {D}_c\mu_{a_2}
      -P_{c}^{\; p}\phi_{pa_2}
      +\rh_{ca_2}  =0.
\end{gather*}

We now consider the map $H:\Lambda^3\TT\to S^2\TT^*$ def\/ined in \eqref{defH}.
As a $\SO(h)$-representation $S^2{\rr^{7}}^*$ decomposes
into the irreducible components $S^2_0{\rr^{7}}^*$ of trace-free
symmetric 2-forms and the space $\rr h$ of multiples of $h$.
The corresponding decomposition on the tractor level is
\begin{gather*}
  S^2\TT^*=S^2_0\TT^*\oplus\rr\mb{h}.
\end{gather*}
Accordingly, $H(\mb{\Phi})$ decomposes into $H(\mb{\Phi})_0$ and
$H(\mb{\Phi})_{tr}$.
Compatibility condition \eqref{compat}\ then means that $H(\mb{\Phi})_0=0$
\ and $H(\mb{\Phi})_{tr}\not=0$.
\begin{lemma}\label{lem0} A parallel tractor $3$-form $\Phi$ satisfies
  $H(\mb{\Phi})_0=0$ if and only if
  \begin{gather*}
    \phi\wedge\phi\wedge\mu=0.
  \end{gather*}
   In particular, $H(\mb{\Phi})_0=0$ whenever $\phi$ is locally decomposable.
\end{lemma}

  \begin{proof}
    $S^2_0\TT^*$ is the tractor-bundle associated to the irreducible representation of $\SO(h)$ on $S^2_0{\rr^{n+2}}^*$. By assumption, $\mb{\Phi}$ is
    $\na^{\La^3\TT}$-parallel. The mapping $\mb{\Phi}\mapsto H(\mb{\Phi})\mapsto H(\mb{\Phi})_0$ is algebraic, and thus naturality of the tractor connection
    implies that $H(\mb{\Phi})_0$ is $\na^{S^2_0\TT^*}$-parallel.
    The section $H(\mb{\Phi})_0\in\Ga(S^2_0\TT^*)$ can thus by Lemma \ref{reproduceparlem}\
    be recovered via the BGG-splitting operator $L_0^{S^2_0\TT^*}$
    from its projection to $\HH_0(S^2_0\TT^*)=\ce[2]$.
    This projection is achieved by inserting twice the top slot
    $\tau_+$ into $H(\mb{\Phi})_0$, but since $\mb{h}(\tau_+,\tau_+)=0$\
    this is the same as evaluating $H(\mb{\Phi})(\tau_+,\tau_+)$.
    Now according to \eqref{defH}\
    \begin{gather*}
      H(\mb{\Phi})(\tau_+,\tau_+)\mb{vol}=(i_{\tau_+}(\mb{\Phi}))\wedge(i_{\tau_+}(\mb{\Phi}))\wedge\mb{\Phi}.
    \end{gather*}
    Using the representation \eqref{splitPhi}\ with
    respect to a metric $g\in[g]$, we have $i_{\tau_+}\mb{\Phi}= \phi-\tau_+\wedge\mu$ and thus
    \begin{gather*}
       H(\mb{\Phi})(\tau_+,\tau_+)\mb{vol}= (\phi-\tau_+\wedge\mu)\wedge(\phi-\tau_+\wedge\mu)\wedge(\tau_-\wedge\phi
       +\ph+\tau_+\wedge\tau_-\wedge\mu+\tau_+\wedge\rh)\\
       \phantom{H(\mb{\Phi})(\tau_+,\tau_+)\mb{vol}}{} =\phi\wedge\phi\wedge\tau_+\wedge\tau_-\wedge\mu
       -\tau_+\wedge\mu\wedge\phi\wedge\tau_-\wedge\phi-\phi\wedge\tau_+\wedge\mu\wedge\tau_-\wedge\phi \\
       \phantom{H(\mb{\Phi})(\tau_+,\tau_+)\mb{vol}}{} =
     3\tau_+\wedge\tau_-\wedge\phi\wedge\phi\wedge\mu.
    \end{gather*}
    This vanishes if and only if $\phi\wedge\phi\wedge\mu=0$.
  \end{proof}

Assume now that $H(\mb{\Phi})_0$ vanishes, i.e. $H(\mb{\Phi})=H(\mb{\Phi})_{tr}=\la \mb{h}$, and since $0=\na^{S^2\TT^*_0}(\la \mb{h})=(d\la)\mb{h}$ we have that $\la\in\rr$\
is a constant.
\begin{lemma}\label{lemtr}
  Suppose that $H(\mb{\Phi})_0=0$. Then, $H(\mb{\Phi})=\la\mb{h}$ for
  a constant $\la\in\rr$, $\la\not=0$,  if and only if
  \begin{gather}\label{condtr}
    \phi\wedge\mu\wedge\rh\not=0.
  \end{gather}
 \end{lemma}

  \begin{proof}
    We check that $\la\not=0$ by inserting $\tau_+$, $\tau_-$ since
$H(\mb{\Phi})(\tau_+,\tau_-)=\la\mb{h}(\tau_+,\tau_-)=\la$:
\begin{gather*}
   H(\mb{\Phi})(\tau_+,\tau_-)\mb{vol}=(i_{\tau_+}\mb{\Phi})\wedge(i_{\tau_-}\mb{\Phi})
  \wedge\mb{\Phi} \\
   \phantom{H(\mb{\Phi})(\tau_+,\tau_-)\mb{vol}}{}
  =(\phi-\tau_+\wedge\mu)\wedge(\tau_-\wedge\mu+\rh)
  \wedge(\tau_-\wedge\phi+\ph+\tau_+\wedge\tau_-\wedge\mu+\tau_+\wedge\rh) \\
  \phantom{H(\mb{\Phi})(\tau_+,\tau_-)\mb{vol}}{} =
  \phi\wedge\tau_-\wedge\mu\wedge\tau_+\wedge\rh+\phi\wedge\rh\wedge\tau_+\wedge\tau_-\wedge\mu-\tau_+\wedge\mu\wedge\rh\wedge\tau_-\wedge\phi \\
  \phantom{H(\mb{\Phi})(\tau_+,\tau_-)\mb{vol}}{}=
  3\tau_+\wedge\tau_-\wedge\phi\wedge\mu\wedge\rh.
\end{gather*}
Thus, $\la\not=0$ if and only if $\phi\wedge\mu\wedge\rh\not=0$. Note
that this f\/ixes the constant $\la$ and $\phi\wedge\mu\wedge\rh$
either vanishes globally or nowhere.
  \end{proof}

We are now ready to prove Theorem~A:
\begin{thma}
  Let $[g]$ be a conformal class of signature $(2,3)$ metrics on $M$.
  Then $[g]$ is induced from a generic rank $2$ distribution $\D\subset TM$\
  if and only if there exists a normal conformal Killing $2$-form $\phi$ which
  is locally decomposable and satisfies the  genericity condition
\begin{gather*}
    \phi\wedge\mu\wedge\rh\not=0.
  \end{gather*}
\end{thma}

\begin{proof}
Suppose $[g]$ is induced from a generic rank $2$ distribution $\D\subset TM$.  Then Proposition~\ref{prop2form}
shows that there is a~parallel tractor $3$-form $\mb{\Phi}\in\Gamma(\La^3\TT)$ that projects to a locally decomposable normal conformal Killing $2$-form $\phi=\Pi_0(\mb{\Phi})$. Furthermore, the discussion in Section~\ref{Phired} shows that $H(\mb{\Phi})$ is a nonzero multiple of the tractor metric $\mathbf{h}.$  This implies $\phi\wedge\mu\wedge\rh\not=0$ by Lemma~\ref{lemtr}.

Conversely, suppose we have a locally decomposable normal conformal Killing $2$-form $\phi\in\Ga(\La^2 T^*M\t\ce[3])$ satisfying the genericity condition. According to Section~\ref{section-cKf}, $\phi$ corresponds to a~parallel tractor $3$-form $\mb{\Phi}$ given by $\eqref{splitPhi}.$ Lemmata~\ref{lem0} and~\ref{lemtr} show that the assumptions on $\phi$ imply $H(\mb{\Phi})=\lambda \mathbf{h},$ for $\lambda\neq 0,$ and thus $\Hol([g])\subset G_2,$ as explained   in Section~\ref{Phired}. By Theorem~\ref{holred}, this means that the conformal structure is canonically associated to a $2$-distribution $\D$.
\end{proof}

\begin{remark}
  It is a well known consequence of the classical Pl\"ucker relations
  (cf.~\cite{eastwood-michor-pluecker}) that a~two form $\ph$ is
  locally decomposable if and only if $\ph\wedge\ph$ vanishes globally.
\end{remark}

\begin{remark}\label{orientability}
Throughout this paper we have assumed orientability of $TM$.
This however, is only a~minor point:
If we leave this assumption and denote by $\mc{O}$ the
$2$-fold covering of $M$ which is the orientation-bundle,
we would obtain a \emph{twisted}\ normal conformal Killing $2$-form
$\phi\in\Ga(\La^2 T^*M\t\ce[3]\t\mc{O})$.
\end{remark}

\section[Decomposition of conformal Killing fields of $\protect{[g]_{\D}}$]{Decomposition of conformal Killing f\/ields of $\boldsymbol{[g]_{\D}}$}\label{section-decomp}

The goal of this section is to prove Theorem B.
We will show that every conformal Killing f\/ield of~$[g]_{\D}$
decomposes into a symmetry of
 the distribution $\D$ and an almost Einstein scale.
The
space of almost Einstein scales $\mb{aEs}([g])$ was def\/ined
in~\eqref{def-aEs}  of Section~\ref{section-aEs}.
Now we discuss symmetries.

Since $\D$ and $[g]$ are equivalently described by Cartan geometries
$(\G,\om)$ resp.\ $(\ti\G,\ti\om)$ we can determine their symmetry algebras
by determining the symmetries
of their corresponding Cartan geometries -- in fact, one can def\/ine them in this way.
For this purpose,
we give the general description~\cite{cap-infinitaut} of the Lie algebra
of inf\/initesimal automorphisms of a parabolic geometry below in Section~\ref{infinitaut}.

Before this Cartan-geometric description, let us discuss the classical
notions. An inf\/initesimal automorphism or symmetry of the distribution
$\D\subset TM$ is a vector f\/ield on $M$ whose Lie derivative preserves
$\D$, i.e.,
\begin{gather}\label{Dsym}
  \mb{sym}(\D)=\{\xi\in\X(M)  :\ \mc{L}_{\xi}\eta=[\xi,\eta]\subset\Ga(\D)\ \forall\, \eta\in\Ga(\D)\}.
\end{gather}

\begin{remark}
  In this text we won't show directly that the symmetries of
the distribution $\mb{sym}(\D)$ def\/ined via~\eqref{Dsym} agree
with the inf\/initesimal automorphisms $\mb{inf}(\om)$ of the corresponding Cartan
geometry discussed below. We just use the fact that associating a regular normal
parabolic geometry of type $(G_2,P)$ to a generic rank $2$-distribution~$\D$
is an equivalence of symmetries.
The explicit form of the splitting from vector f\/ields on $M$ into the
adjoint tractor bundle  relating the classical and the Cartan-viewpoint
is only given in the conformal case, since there we will later need the
explicit formula.
\end{remark}

\subsection{Conformal Killing f\/ields}
The symmetries of the associated conformal structure $[g]=[g]_{\D}$
are the conformal\ Killing f\/ields
\begin{gather*}
  \mb{cKf}([g])=\{\xi\in\X(M)  :\ \mc{L}_{\xi} g=e^{2f}g\ \mr{for\ some}\ g\in[g]\ \mr{and}\ f\in\cinf(M)\}.
\end{gather*}
Since $\mc{L}_{\xi}g$ decomposes into
a multiple of $g$ and a trace-free part one can
equivalently demand that $\mc{L}_{\xi}g$ is pure trace.
Now, with $D$ the Levi-Civita connection of $g\in [g]$,  $\mc{L}_{\xi}g$ being pure trace is equivalent to
\begin{gather*}
  D_{(c}\xi_{a)_0}=D_{(c} g_{a)_0p}\xi^p=0;
\end{gather*}
i.e., the symmetric, trace-free part of $D_{c}\xi_a$ vanishes.

As an equation on $1$-forms of conformal weight $2$ this
is in fact described by the f\/irst BGG-operator of $\La^2 \TT$:
By \eqref{L0form}, the splitting operator
\[L_0^{\La^2\TT}: \ \X(M)=\ce_a[2]\goesto \La^2\TT=\ti\A M\]  is given by
\begin{gather}\label{L0La2}
  \si_a\in\ce_a[2]\mapsto
  \begin{pmatrix}
    \begin{pmatrix}
                -\frac{1}{10}D^pD_p\si_{a}+\frac{1}{10}D^pD_{a}\si_{p} -\frac{1}{5}J\si_{a}
\\
+\frac{2}{5}P^p_{\; a}\si_{p}+\frac{1}{25}D_{a}D^p\si_{p}
    \end{pmatrix}\\
    D_{[c}\si_{a]} \ | \ -\frac{1}{5}D^p\si_p   \\
    \si_a
  \end{pmatrix}.
\end{gather}
Now the f\/irst BGG-operator $\Th^{\La^2\TT}_0$ of $\La^2\TT$
def\/ined by the composition
\begin{gather*}
  \Th^{\La^2\TT}_0 : \ \X(M)=\ce_a[2]\goesto\ce_{(ab)_0}[2], \qquad
  \Th^{\La^2\TT}_0 =\Pi_1\circ\na^{\La^2\TT}\circ L_0^{\La^2{\TT}}
\end{gather*}
is seen by direct calculation employing \eqref{tracon}
to be
\begin{gather*}
  \xi^a\mapsto D_{(c}\xi_{a)_0}
\end{gather*}
for $\xi\in\X(M)$; i.e., $\Th_0^{\La^2\TT}$\
is the conformally invariant operator governing
conformal Killing f\/ields.

We now proceed to prove a technical lemma to be used in the proof of Theorem \ref{thm-decomp}\ below.
It is a general fact  for  parabolic geometries satisfying a certain homological condition that parallel sections of the adjoint tractor bundle insert trivially into curvature~\cite{cap-infinitaut}, but it is easy to see this directly for conformal Killing f\/ields, where this has f\/irst been observed in~\cite{gover-lapl_einstein}.
We only sketch a simple proof for conformal structures
of dimension $\geq 4$, which is all we need:
\begin{lemma}\label{lem-trivinsertion}
  Let $s\in\Ga(\ti\A M)$ such that $\na^{\ti\A}s=0$.
  Then $\ti K(\Pi(s),\cdot)=0$.
\end{lemma}

  \begin{proof}
    Since $\na^{\ti A}s=0$ one has $R_{c_1c_2}s=\na^{\ti A}_{[c_1}\na^{\ti\A}_{c_2]}s=0$, with $R\in\Om^2(M,\gl(\ti\A M))$ the curvature of~$\na^{\ti\A}$.
    But since $\na^{\ti\A}$ is the induced tractor connection of
    the Cartan connection form~$\ti\om$, one has that
      $Rs$ is the algebraic action of $\ti K\in\Om^2(M,\ti\A M)=\Om^2(M,\so(\TT))$ on~$s$; thus, $\ti K$ annihilates~$s$.
    Via the projection $\Pi:\Ga(\La^2\mb{T})\goesto\ce_a[2]=\X(M)$,
    $s$ projects to a conformal Killing f\/ield~\mbox{$\si^a\in\X(M)$},
    and from the explicit formula~\eqref{formulaK} one obtains that
    the Weyl curvature~$C$
    annihilates~$\si$. Employing the symmetries of~$C$,
    one then immediately has that this is equivalent to trivial
    insertion of~$\si$. Then $\ti K(\si,\cdot)=0$
    follows from $(n-2)A_{abc}=D^pC_{pabc}$.
  \end{proof}

\subsection{Inf\/initesimal automorphisms of parabolic geometries}\label{infinitaut}

In this subsection we will use the description~\cite{cap-infinitaut}
of the symmetry Lie algebra
of a parabolic geometry $(\G,\om)$ of type $(G,P)$. This will be applied
for $(\G,\om)$ being, as above, the geometry of type $(G_2,P)$\
describing the generic rank 2 distribution $\D$ and for
the conformal geometry encoded in the Cartan geometry $(\ti\G,\ti\om)$\
of type $(\SO(h),\ti P)$. For
more details see \cite{cap-infinitaut}\ or \cite{mrh-hcg}.

Since $\G\goesto M$ is a $P$-principal bundle over $M$ and the geometric
structure is encoded in the Cartan connection form $\om\in\Om^1(\G,\g)$,
one def\/ines an automorphism of $(\G,\om)$ as a
$P$-equivariant dif\/feomorphism $\Psi$ of $\G$ preserving~$\om$, i.e., $ \Psi^*\om=\om$.
It is well known that the automorphism group of a Cartan geometry is
always a Lie group (see for instance \cite{cap-slovak-par}). The space of \emph{infinitesimal automorphisms} is given by the  $P$-invariant vector f\/ields on $\G$ preserving~$\om$, i.e.,
 \begin{gather*}
   \mb{inf}(\G,\om)=\{\xi\in\X(M)^P : \ \mc{L}_{\xi}\om=0\}.
 \end{gather*}
The Lie algebra of the automorphism group of $(\G,\om)$ then consists of those  vector f\/ields in $\mb{inf}(\G,\om)$ that are complete.

Now $\om:T\G\goesto \g$ is a $P$-equivariant trivialization of $T\G$,
and thus, for a $\xi\in\mb{inf}(\G,\om)$ the function $f=\om\circ\xi:\G\goesto\g$\
is $P$-equivariant. The function $f$ therefore def\/ines a section of
the adjoint tractor bundle $\A M=\G\times_P \g$.

We have an explicit formula for the tractor connection $\na^{\A}$ on $\A M$, see
\eqref{deftracon}: Let $s\in\Ga(\A M)$ be the section corresponding to the $P$-equivariant function $f\in\cinf(\G,\g)^P$.
To compute $\na^{\A}_{\eta}s$ for $\eta\in\X(M)$\
we take a $P$-invariant lift $\eta'\in\X(\G)$ of $\eta$. Then $\na^{\A}_{\eta}s$ corresponds to the $P$-equivariant map
\begin{gather*}
  u\mapsto \eta'_u\cdot\om(\xi)+[\om_u(\eta'),\om_u(\xi)].
\end{gather*}
For the next lemma, f\/irst note that the natural projection
\begin{gather*}
  \Pi: \ \A M=\G\times_P\g\goesto\G\times_P\g/\p= TM
\end{gather*}
projects $s\in\Ga(\A M)$  to a vector f\/ield $\Pi(s)\in\X(M)$, which is in fact just the projection of the $P$-invariant vector f\/ield $\xi\in\X(\G)^P$.
Thus one can insert $\Pi(s)\in\X(M)$ into the curvature form $K\in\Om^2(M,\A M)$.
\begin{proposition}[\cite{cap-infinitaut}]\label{prop-infaut}
  Let $s\in\Ga(\A M)$ be the adjoint tractor corresponding to
  the $P$-invariant vector field $\xi\in \X(\G)^P$.
  Then
  \begin{gather*}
    \mc{L}_{\xi}\om=0\ \mr{if\/f}\ \na^{\A}_{\eta}s+K(\Pi(s),\eta)=0\ \forall\eta\in\X(M),
  \end{gather*}
  i.e., $\xi\in\mb{inf}(\G,\om)$ if and only if  the corresponding
  adjoint tractor $s\in\Ga(\A M)$ is parallel with respect to the
  connection
  \begin{gather*}
    \hat\na^{\A}s=\na^{\A}s+K(\Pi(s),\cdot).
  \end{gather*}
\end{proposition}

\subsection{Decomposition of the conformal adjoint tractor bundle.}
Let $\A M:=\G\times_P\g$ be the adjoint tractor bundle of a
generic 2-distribution $\D$ and let $\ti\A M:=\ti\G\times_{\ti P}\so(h)$ be the
adjoint tractor bundle of the associated conformal structure.
The tractor connection on $\A M$ will be denoted by $\na^{\A}$
and the one on $\ti\A M$ by $\na^{\ti\A}$.

Recall from Section~\ref{somealgebra} that as a $G_2$-representation,
  $\so(h)=\rr^7\oplus\g$,
  i.e., $\so(h)$ decomposes into the direct sum of the
  standard representation of $G_2\subset\SO(h)$ on $\rr^7$ and
  the adjoint representation $\Ad:G_2\goesto\GL(\g)$.
This decomposition was realized by the exact sequence~\eqref{sequence-Phi} of Section~\ref{somealgebra} and its splitting~\eqref{split2}.
On the level of associated bundles it yields a decomposition
of the conformal adjoint tractor bundle:
\begin{gather}\label{adj-decomp}
  \ti\A M=\ti\G\times_{\ti P}\so(h)=\G\times_P\so(h)=(\G\times_P\rr^7)
  \oplus(\G\times_P\g)=\TT\oplus\A M.
\end{gather}
I.e., the conformal adjoint tractor bundle $\ti\A M$ is the direct sum
of the standard tractor bundle~$\TT$ and the adjoint tractor bundle $\A M$
of the generic distribution.

Let us check that we can also decompose the tractor connection into
\begin{gather}\label{tra-decomp}
  \na^{\ti\A}=\na^{\TT}\oplus\na^{\A}.
\end{gather}
Take a vector f\/ield $\xi\in\X(M)$ and its horizontal lift $\xi'$
to a vector f\/ield on the extended bundle $\G'=\G\times_P G_2$.
Then, for $s\in\Ga(\ti\A M)$, the tractor derivative $\na_{\xi}s$
is def\/ined by dif\/ferentiating the $G_2$-equivariant function
$f:\G\goesto\so(h)$ corresponding to $s$ in direction $\xi'$.
But taking this derivative evidently commutes with the algebraic
projections of $f$ to its components $f_{\TT}:\G\goesto\rr^7$ and $f_{\A M}:\G\goesto\g$; thus \eqref{tra-decomp}\ holds.

To prove the decomposition of conformal Killing f\/ields we need the following theorem about the deformed connections
  \begin{gather*}
    \hat\na^{\A}s=\na^{\A}s+K(\Pi(s),\cdot)
  \end{gather*}
  and
    \begin{gather*}
    \hat\na^{\ti\A}s=\na^{\ti\A}s+\ti K(\Pi(s),\cdot)
  \end{gather*}
  whose parallel sections describe inf\/initesimal automorphisms.
 Recall that according to Lemma~\ref{lem-curvaturerelations}
  we have $\ti K=K \in\Om^2(M,\A M)$.

\begin{theorem}\label{thm-decomp}
  Let $s\in\Ga(\ti\A M)$ be a section of the conformal adjoint tractor bundle which splits according to decomposition
  \eqref{adj-decomp}\ into $s_1\in\Ga(\TT)$ and $s_2\in\Ga(\A M)$.
  Then $s$ is parallel with respect to $\hat\na^{\ti\A}$ if and only if
  $s_1$ is $\na^{\TT}$-parallel and $s_2$ is $\hat\na^{\A}$-parallel.
\end{theorem}

\begin{proof}
   Let $s_1\in\Ga(\TT)$ and $s_2\in\Ga(\A M)$ be $\na^{\TT}$- resp.\
  $\hat\na^{\A}$-parallel sections.
  Since the restriction of $\na^{\ti\A}$ to $\Ga(\TT)\subset\Ga(\ti\A M)$\
  is just $\na^{\TT}$, the section $s_1$ includes as a $\na^{\ti\A}$-parallel section into $\Ga(\ti\A M)$.
  Lemma \ref{lem-trivinsertion}\ shows that that $K(\Pi(s_1),\cdot)=0$, and thus
  also $\hat\na^{\ti\A}s_1=0$. For $s_2$ we have $\hat\na^{\ti\A}s_2=0$ by Lemma
\ref{lem-curvaturerelations},
  and thus $s:=s_1+s_2$  satisf\/ies $\hat\na^{\ti\A}s=0$.

  Conversely, we take  $s\in\Ga(\ti\A M)$ with $\hat\na^{\ti\A}s=0$ and
  decompose $s=s_1\oplus s_2\in\Ga(\TT)\oplus\Ga(\A M)$ according to \eqref{adj-decomp}.
  Since $K$ has values in $\A M$\
  we have that $s_1\in\Ga(\TT)$ is parallel with respect
  to the standard tractor connection $\na^{\TT}$ by \eqref{tra-decomp}.
  We still need to show that $s_2$ is parallel with respect to
  $\hat\na^{\A}$, while so far we know that
  \begin{gather*}
    \na^{\A}s_2+K(\Pi(s_1),\cdot)+K(\Pi(s_2),\cdot)
  \end{gather*}
  vanishes.
  But since $s_1$ is parallel as a section of $\ti\A M$ with
  respect to the usual adjoint tractor connection $\na^{\ti \A}$ according
 to \eqref{tra-decomp}, we can
  again apply Lemma \ref{lem-trivinsertion}, which
  tells us that $s_1$ inserts trivially into the curvature $\ti K=K$.
  Thus also $\hat\na^{\A}s_2=0$.
\end{proof}

Proposition \ref{prop-infaut}  provides an identif\/ication
  of $\mb{inf}(\ti\G,\ti\om)$ with parallel sections of $\ti\A M$\
  with respect to $\hat\na^{\ti \A}$
  and an identif\/ication of $\mb{inf}(\G,\om)$ with parallel sections
  of $\A M$ with respect to the connection $\hat\na^{\A}$.
We can thus translate Theorem~\ref{thm-decomp}
into a decomposition of conformal Killing f\/ields:

\begin{thmb}
      Let $[g]$ be conformal class of signature $(2,3)$-metrics
      on $M$ and $\phi\in\ce_{[ab]}[3]$ a locally decomposable
      normal conformal Killing two form that
      satisfies the genericity condition~\eqref{condtr}.
      By Theorem~A, there is a generic distribution $\D$
      with $[g]=[g]_{\D}$.

      Then every conformal Killing field decomposes into a symmetry of
    the distribution $\D$ and an almost Einstein scale:
    \begin{gather}\label{decomp}
      \mb{cKf}([g])=\mb{sym}(\D)\oplus\mb{aEs}([g]).
    \end{gather}

    The mapping that associates to a conformal Killing field $\xi\in\X(M)$\
    its
    almost Einstein-scale part with respect to the decomposition~\eqref{decomp}
    is given by
\begin{gather}\label{einsteinpart}
      \xi_a\mapsto \phi_{pq}(D\xi)^{pq}-\tfrac{1}{2}\xi^p D^q \phi_{pq},
    \end{gather}
    where $D$ is the Levi Civita connection of an arbitrary metric
    $g$ in the conformal class.

    The mapping that associates to an almost Einstein scale $\si\in\ce[1]$\
    a conformal Killing field is given by
\begin{gather}\label{einsteinfield}
       \si\mapsto \phi_{ap}D^p\si-\tfrac{1}{4}\si D^p\phi_{pa}.
    \end{gather}
\end{thmb}

    \begin{proof}
By Proposition \ref{prop-infaut} conformal Killing f\/ields
of $[g]$ are in 1:1-correspondence with $\hat\na^{\ti \A}$-parallel
sections of $\ti\A M$. By the theorem above every such
section decomposes into a parallel standard tractor
in $\Ga(\TT)$ and a  $\hat\na^{\A}$-parallel section of $\A M$.
By Proposition \ref{prop-aEs}\ and again Proposition \ref{prop-infaut}, now
for $\A M$, this yields the decomposition \eqref{decomp}.

It is now straightforward to make this decomposition explicit in
terms of the normal conformal Killing $2$-form of Theorem~A
encoding the generic distribution $\D$.

To map an almost Einstein scale $\si\in\ce[1]$ to a conformal
Killing f\/ield we use the splitting operator
$L_0^{\TT}:\ce[1]\goesto \Ga(\TT)$ given in \eqref{splitStd},
contract this section into the characterizing $\na^{\La^3 \TT}$-parallel $3$-form
$\mb{\Phi}\in\Ga(\La^3\TT)$ given by~\eqref{splitPhi}
via the tractor metric $\mb{h}$,
and project the resulting section of $\La^2\TT=\ti\A M$ down
to $\X(M)$. This yields~\eqref{einsteinfield}.

To project a conformal Killing f\/ield $\xi\in\X(M)$\
to its almost Einstein scale-part we proceed similarly: we map it to $\La^2\TT$ via \eqref{L0La2}, contract it into $\mb{\Phi}\in\Ga(\La^3\TT)$ and project the resulting
standard tractor to $\ce[2]$. This gives \eqref{einsteinpart}.
    \end{proof}

To be precise, one has to use a constant scalar multiple of $\phi\in\ce_{[ab]}[3]$\
such that the composition of \eqref{einsteinpart}\ with \eqref{einsteinfield}\ is
the identity.

\begin{remark}
Mapping \eqref{einsteinfield}\ actually works more generally: in the presence
of an almost Einstein scale it was shown in \cite[Corollary~5.2]{gover-silhan-2006} that one can associate to every conformal Killing $2$-form, not only
to normal ones, a conformal Killing f\/ield.
\end{remark}

\subsection*{Acknowledgements}
We cannot underestimate the value of discussions with Andreas \v{C}ap
on various technical procedures used in this paper. The concept of orbit types
of parallel tractors was introduced to the f\/irst author by Felipe
Leitner, who moreover suggested to check for simplicity of the
underlying $2$-form. We thank the referees for their careful reading
and various valuable suggestions for improvements.

The f\/irst author was supported by the IK I008-N funded by the University of Vienna. The second author was supported by project  P 19500-N13 of the ``Fonds zur F\"orderung der wissenschaftlichen Forschung'' (FWF).

\pdfbookmark[1]{References}{ref}

\LastPageEnding

\end{document}